\documentclass[12pt,a4paper]{amsart}

\usepackage{euscript,amsfonts,amssymb,amsmath,amscd,amsthm,enumerate,hyperref}

\usepackage{comment}

\usepackage{mathrsfs}

\usepackage{graphicx}

\usepackage{color}

\usepackage[margin=1.1in]{geometry}

\usepackage{bbm}

\newtheorem{theorem}{Theorem}[section]
\newtheorem*{theorem*}{Theorem}
\newtheorem{lemma}[theorem]{Lemma}
\newtheorem{definition}[theorem]{Definition}
\newtheorem{proposition}[theorem]{Proposition}

\theoremstyle{remark}
\newtheorem{rmk}[theorem]{Remark}

\newcommand{\E}{\mathbb E}
\newcommand{\EE}{\mathbb E}

\newcommand{\pp}{\mathbb{P}}

\newcommand{\rr}{\mathbb{R}}

\newcommand{\nn}{\mathbb{N}}
\newcommand{\ep}{\hfill \ensuremath{\Box}}
\newcommand{\eq}{\begin{equation}}
\newcommand{\en}{\end{equation}}

\newcommand{\bone}{\mathbf{1}}
\newcommand{\PP}{\mathbb{P}}
\newcommand{\RR}{\mathbb{R}}

\numberwithin{equation}{section}

\title[Particle systems]{Particle systems with singular interaction through hitting times: application in systemic risk modeling} 

\author{Sergey Nadtochiy}
\address{Department of Mathematics, University of Michigan, Ann Arbor, MI, USA}
\email{sergeyn@umich.edu}
\thanks{S. Nadtochiy is partially supported by the NSF grant DMS-1411824.}

\author{Mykhaylo Shkolnikov}
\address{ORFE Department, Princeton University, Princeton, NJ, USA}
\email{mshkolni@gmail.com}
\thanks{M. Shkolnikov is partially supported by the NSF grant DMS-1506290.}

\begin{document}

\begin{abstract}
We propose an interacting particle system to model the evolution of a system of banks with mutual exposures. In this model, a bank defaults when its normalized asset value hits a lower threshold, and its default causes instantaneous losses to other banks, possibly triggering a cascade of defaults. The strength of this interaction is determined by the level of the so-called \emph{non-core exposure}.
We show that, when the size of the system becomes large, the cumulative loss process of a bank resulting from the defaults of other banks exhibits discontinuities. These discontinuities are naturally interpreted as \emph{systemic events}, and we characterize them explicitly in terms of the level of non-core exposure and the fraction of banks that are ``about to default".
The main mathematical challenges of our work stem from the very singular nature of the interaction between the particles, which is inherited by the limiting system. A similar particle system is analyzed in \cite{Delarue1} and \cite{Delarue2}, and we build on and extend their results. In particular, we characterize the large-population limit of the system and analyze the jump times, the regularity between jumps, and the local uniqueness of the limiting process.
\end{abstract}

\maketitle

\section{Introduction}

Consider an interconnected system whose components might fail, such as a banking system in which banks may default. The existing approaches to quantitative modeling of such systems, roughly, fall into the following two categories: (1) network models, and (2) particle systems with mean-field interaction. 
Models of the first category are considered, e.g., in \cite{May}, \cite{Gai}, \cite{Lorenz}.
These models are able to capture the current characteristics of the system with high precision and to predict the effects of immediate external shocks.
However, to obtain analytical results on the risks associated with a given network (for example, on the probability of a default cascade of a certain size due to a specific external shock) a limit, as the size of the system goes to infinity, needs to be taken. The results are, then, expressed in terms of the average values of the members' characteristics. 
In the models of the second category, it is assumed from the very beginning that there are only a few characteristics by which any two particles (e.g. banks) can differ and that the interaction between them is either of a mean-field type, or is given implicitly through a correlation with a common factor. In contrast to the network models, the particle systems are dynamic and allow to investigate how the system evolves over time. Under suitable assumptions, it is possible to derive analytic formulas for the risk of future failures in such systems, in addition to describing the immediate (i.e. static) risk embedded in the system at a given point in time. Analysis of this kind is often carried out in the context of losses due to defaults in a large portfolio, see, e.g., \cite{DuffieDembo}, \cite{Giesecke}, \cite{Hambly}, \cite{Fouque}, \cite{Horst}, \cite{DaiPra}, and the references therein.

\medskip

In the present paper, we follow the mean-field approach in modeling the dynamics of an interconnected system of banks. At the same time, we use an explicit (i.e. structural) mechanism of default contagion (i.e. interaction between the particles), which, in particular, differentiates our setting from the models based on interacting default intensities (cf., \cite{Giesecke}, \cite{DaiPra}). The goal of our paper is to provide a method for estimating the proximity of a \emph{systemic failure} (i.e. for quantifying the \emph{systemic risk}), which would allow a regulator to intervene ahead of time.
We understand the systemic failure as the occurrence of a ``significantly large" default cascade. It has been documented in various studies (cf., \cite{Gai}, \cite{Lorenz}, \cite{DaiPra}) that an interconnected banking system transitions between two regimes: the well-behaved regime, in which the system spends most of its time, and during which the default cascades are very small or do not appear at all, and the systemic crisis regime, which occurs rarely, and which is characterized by large groups of banks defaulting in a short period of time. Even though the presence of such \emph{phase transitions} is well known, it is often difficult to define precisely what constitutes a significantly large cascade. In our setting the times of such cascades are captured by the discontinuity points of the cumulative loss process in a limiting system, thus, providing a natural endogenous definition of a systemic failure.

\medskip

In addition, we provide an explicit connection between the occurrence of systemic events and the internal characteristics of the banking system, which, in principle, can be observed and controlled by a regulator. More specifically, we describe the time of systemic failure in terms of the level of \emph{mutual exposure} of the banks and the fraction of banks in the immediate danger of default. The level of mutual exposure measures the interconnectedness of the system, allowing the firms to lend to and borrow from each other. Such lending and borrowing may decrease the individual risk of a bank, but it also provides channels for spreading losses from individual defaults across the rest of the system (i.e. for creating default contagion). This dual role of mutual exposure has been analyzed, e.g., in \cite{BattistonStiglitz}, \cite{Gai}, \cite{Lorenz}, \cite{Glasserman}.
The mutual exposure of the banks is also known as the \emph{non-core exposure} (in contrast to the core exposure, which measures how much the banks lend to the real economy), and its effects on the occurrence of systemic crises is analyzed, e.g., in \cite{Shin}, \cite{DuffieZhu}.
As mentioned above, the level of non-core exposure can be controlled by a financial regulator, and such a control was implemented by the government of South Korea in the aftermath of the financial crisis of 2008 (cf. \cite{Shin}, \cite{IMF}).
The joint effect of the interconnectedness of the system (measured by the level of non-core exposure, in the present case) and the fraction of members in the immediate danger of failure on the occurrence of a large failure cascade has been also investigated in \cite{BattistonStiglitz}, \cite{Lorenz}, \cite{Watts}.


\medskip

The mathematical model considered in this paper is based on a system of particles with singular interaction through hitting times. Systems of this type are considered, e.g., in \cite{Delarue1} and \cite{Delarue2} in connection to a problem from neuroscience, and we use and extend some of the ideas developed therein to establish our results (for alternative models, with smooth interaction, we refer to \cite{Talay} and \cite{Giesecke}). In particular, the convergence of the $N$-particle system as $N\to\infty$ (Theorem \ref{main2}) is based on appropriate adaptations of the methods presented in \cite{Delarue2}. The main original contribution of our work is in the analysis of the limiting process. Namely, Theorem \ref{main2.2} provides a sufficient condition for a jump in the cumulative loss process, Theorem \ref{thm:Holder.loss} establishes regularity of the limiting process, while Theorem \ref{main1} proves its local uniqueness. Theorem \ref{main1} amounts to proving the uniqueness of the solution to a non-local non-linear Cauchy-Dirichlet problem. It is worth mentioning that the limiting process of a similar particle system has been analyzed in \cite{Delarue1} (see also \cite{Perthame}, \cite{Perthame2}, where the focus is on documenting the possibility of a jump in the limiting process and on describing stationary solutions when no jumps occur). The paper \cite{Delarue1} establishes regularity of the limiting process and proves its uniqueness using a similar Cauchy-Dirichlet problem. However, the main results of \cite{Delarue1} require additional assumptions on the strength of interaction in the system, which, in particular, rule out the possibility of a jump in the limiting process.
As such jumps have a natural practical interpretation (e.g. as systemic crises, in the application proposed herein), we, specifically, focus on the cases where such jumps may occur. As a consequence, the limiting process and the solution to the associated Cauchy-Dirichlet problem, herein, do not possess as much regularity as in \cite{Delarue1}, which, naturally, complicates the analysis.

\section{Main results}

Consider $N$ banks and write $X_t^1,\,X_t^2,\,\ldots,\,X_t^N$ for their total asset values at a time $t\ge0$, discounted according to the (possibly stochastic) growth rate of the overall banking system. A bank $i$ defaults when its total asset value $X^i$ drops below a barrier $\underline{x}^i>0$. Since each asset value process can be normalized by the corresponding barrier, we may assume without loss of generality that $\underline{x}^1=\underline{x}^2=\cdots=\underline{x}^N=1$. In the \textit{absence of defaults} we let the asset value processes $X_t^1,\,X_t^2,\,\ldots,\,X_t^N$ follow the stochastic differential equations (SDEs) 
\begin{equation}\label{dynamics_no_defaults}
\begin{split}
\mathrm{d}X_t^i = X^i_t\,(\alpha + \sigma^2/2)\,\mathrm{d}t + X^i_t\,\sigma\,\mathrm{d}B_t^i,\quad i=1,\,2,\,\ldots,\,N,
\end{split}
\end{equation}
where $\alpha\in\RR$ and $\sigma>0$ are constants, and $B^1,\,B^2,\,\ldots,\,B^N$ are independent standard Brownian motions. Here $\alpha$ stands for the return associated with the \textit{traditional investments} of a bank (that is, investments in companies outside of the banking system), and $\sigma$ is the volatility coefficient of such investments. 

\medskip

Now, suppose that the asset value process of a bank $i$ hits the default barrier $\underline{x}^i=1$ at a time $t$, leading to the default of bank $i$. As a result of the default, the asset values of other banks drop and may immediately cause further defaults, and so on. When the default of a bank causes immediate defaults of other banks, we speak of a \textit{default cascade}. For the sake of tractability we assume that, if $k$ banks default at time $t$, then the value of each remaining bank is reduced by the factor 
$$
\left(1-\frac{k}{S_{t-}}\right)^{-C},
$$
where $S_t$ is the number of banks that have survived up to and including time $t$, and $C\in [0,1)$ is a fixed constant. The value of $C$ represents the level of non-core exposure in the banking system. Notice that $\left(1-k/S_{t-}\right)^{C}\approx 1 - C k/S_{t-}$, if $k/S_{t-}$ is small. In such a case, the proposed loss function represents the losses from default contagion in a banking system in which every bank, in total, borrows from all other banks the fraction $C$ of the average bank's value, with the sizes of individual loans distributed proportionally to the other banks' values.
After the default event, the asset value processes of the surviving banks continue to follow the dynamics of \eqref{dynamics_no_defaults} until one of them hits $1$, and so on. 

\medskip

The informal description of the processes $X^1,\,X^2,\,\ldots,\,X^N$ in the previous two paragraphs can be formalized as follows. We fix a time horizon $T>0$, let  
\begin{equation}\label{eq.Y.as.logX.def}
Y^i:= \log X^i,
\quad \tau^i:=\inf\,\{t\in[0,T]:\,Y^i_t\le 0\},\quad i=1,\,2,\,\ldots,\,N
\end{equation}
be the \textit{logarithmic asset values} and the \textit{default times} of the banks, respectively, and denote by 
\begin{equation}\label{eq.V.asMassAlive.def}
S_t:=\sum_{i=1}^N \mathbf{1}_{\{\tau^i>t\}},\quad t\in[0,T]
\end{equation}
the \textit{size of the banking system}. In addition, for any fixed $t\in[0,T]$ we consider the order statistics
\begin{equation}
Y^{(1)}_{t-}\le Y^{(2)}_{t-}\le\cdots\le Y^{(S_{t-})}_{t-}
\end{equation}
of the vector $(Y_{t-}^i:\,\tau^i\ge t)$. Then, the number of defaults at time $t\in[0,T]$ is defined by 
\begin{equation}\label{eq.V.asCascadeSize.def}
D_t=\left(\inf\,\bigg\{k=1,\,2,\,\ldots,\,S_{t-}:\;Y^{(k)}_{t-} + C\log\bigg(1-\frac{k-1}{S_{t-}}\bigg) > 0\bigg\}-1\right)\wedge S_{t-},
\end{equation} 
with the convention $\inf\emptyset=\infty$. Finally, each of the processes $Y^1,\,Y^2,\,\ldots,\,Y^N$ satisfies
\begin{equation}
Y^i_t=\widetilde{Y}^i_t\,\mathbf{1}_{\{\widetilde{Y}^i_s>0,\,s\in[0,t)\}}+\widetilde{Y}^i_{\tau^i}\,\mathbf{1}_{\{t>\tau^i\}},
\end{equation}
where
\begin{equation}\label{dynamics_w_defaults.2.1}
\widetilde{Y}^i_t=Y^i_0+\alpha\,t+\sigma\,B^i_t
+\big(1\wedge(\widetilde{Y}^i_t+1)^+\big) \sum_{u\le t:\,D_u>0} C\log\bigg(1-\frac{D_u}{S_{u-}}\bigg)
\end{equation}
for $t\in[0,\tau^0)\cap[0,T]$,
\begin{equation}\label{dynamics_w_defaults.2.2}
\widetilde{Y}^i_t= (-1)\wedge \widetilde{Y}^i_{\tau^0-} + \alpha\,(t-\tau^0)+\sigma\,(B^i_t - B^i_{\tau^0})
\end{equation}
for $t\in[\tau^0,\infty)\cap[0,T]$, and $\tau^0:=\max_{1\le j\le N} \tau^j$. 

\medskip

It is easy to see that the fixed-point equations defining the auxiliary processes $\widetilde{Y}^1,\,\widetilde{Y}^2,\,\ldots,\,\widetilde{Y}^N$, in (\ref{dynamics_w_defaults.2.1}), are uniquely solvable. Thus, the paths of the processes $\widetilde{Y}^1,\,\widetilde{Y}^2,\,\ldots,\,\widetilde{Y}^N$ and $Y^1,\,Y^2,\,\ldots,\,Y^N$ can be constructed sequentially on the time intervals from $0$ to the first default time, from the first default time to the next default time, and so on. The truncation factors $1\wedge(\widetilde{Y}^i_t+1)^+$, $i=1,\,2,\,\ldots,\,N$, included for a purely technical reason, ensure that each process $\widetilde{Y}^i$ does not jump below $-1$. Nevertheless, these factors are constantly equal to $1$ on $[0,\tau^i)$, for $i=1,\,2,\,\ldots,\,N$, respectively, and therefore have no effect on the pre-default paths $Y^i_t$, $t\in[0,\tau^i)$, $i=1,\,2,\,\ldots,\,N$. In addition, we extend the path of each $\widetilde{Y}^i$ to the interval $[0,T+1]$ continuously, as follows:
\begin{equation*}
\widetilde{Y}^i_t = \widetilde{Y}^i_T + \alpha\,(t-T) + \sigma\,(B^i_t - B^i_T),\quad t\in(T,T+1].
\end{equation*}
The reason we continue the path of $\widetilde{Y}^i$ beyond $\tau^i\wedge T$, rather than stop the process at this time, is that the paths need to be sufficiently ``noisy" in order to establish the desired convergence result. 

\medskip

We sometimes refer to the vector of processes $(\widetilde{Y}^1,\widetilde{Y}^2,\ldots,\widetilde{Y}^N)$ as the \textit{finite-particle system}. It is worth noting that 
\begin{equation}
\sum_{u\le t:\,D_u>0} C\log\bigg(1-\frac{D_u}{S_{u-}}\bigg)
=\sum_{u\le t:\,D_u>0} C\,(\log S_u-\log S_{u-})
=C\log\frac{S_t}{N},\;\;t\in[0,\tau^0\wedge T),
\end{equation}
hence,
\begin{equation}\label{dynamics.YN1}
\widetilde{Y}^i_t=Y^i_0+\alpha\,t+\sigma\,B^i_t+\big(1\wedge(\widetilde{Y}^i_t+1)^+\big)\,C\log\bigg(\frac{1}{N}\,\sum_{j=1}^N \mathbf{1}_{\{\tau^j>t\}}\bigg),\;\;t\in[0,\tau^0\wedge T).
\end{equation}
However, the strong solution of \eqref{dynamics.YN1} is not unique, because the default cascades are not uniquely determined by \eqref{dynamics.YN1} alone. 


\medskip

Being interested in the emergence of large (``systemic'') losses due to default cascades, we study the large $N$ asymptotics of the banking system by means of the empirical measures 
\begin{equation}
\mu^{N}:=\frac{1}{N}\,\sum_{i=1}^N \delta_{Y^i}\quad\text{and}\quad\widetilde{\mu}^{N}:=\frac{1}{N}\,\sum_{i=1}^N \delta_{\widetilde{Y}^i}. 
\end{equation}
We view $\mu^{N}$ and $\widetilde{\mu}^{N}$ as random probability measures on the spaces $D([0,T])$ and $D([0,T+1])$, respectively, consisting of real-valued c\`adl\`ag paths. The latter are endowed with the Skorokhod M1 topology (see, e.g., \cite{Delarue2}, \cite{Skorohod}, \cite{JacodShiryaev}, \cite{Whitt}, for a detailed discussion of the M1 topology). The limiting object associated with the sequence $\widetilde{\mu}^N$, $N\in\nn$ turns out to be given by the following definition (see Theorem \ref{main2} below).

\begin{definition}\label{def_phys_sol.new}
A real-valued c\`adl\`ag process $\overline{Y}_t$, $t\in[0,T]$ is called a \emph{physical solution} if, with
\begin{equation}
\overline{\tau}^0:=\inf\big\{t\in[0,T]:\,\PP\big(\inf_{s\in[0,t]}\overline{Y}_s\leq 0\big)=1\big\},
\end{equation}
it holds 
\begin{equation}\label{eq:phys_sol.1}
\begin{split}
&\overline{Y}_t=\overline{Y}_{0} + \alpha\,t+\sigma\,\overline{B}_t+\big(1\wedge(\overline{Y}_t+1)^+\big)\,\Lambda_t,\;\;t\in[0,\overline{\tau}^0)\cap[0,T], \\  
& \Lambda_t=C\log\,\PP(\overline{\tau}>t),\;\;t\in[0,\overline{\tau}^0)\cap[0,T],\\ 
& \overline{\tau}=\inf\,\{t\in[0,T] : \,\overline{Y}_t\le 0\},
\end{split}
\end{equation}
where $\overline{B}$ is a standard Brownian motion independent of $\,\overline{Y}_0\!>\!0$; and, whenever $\Lambda_t\!\neq\!\Lambda_{t-}$ for some $t\in[0,\overline{\tau}^0)\cap[0,T]$, we have $\Lambda_{t-} - \Lambda_{t} \leq F_t(\overline{D}_t)$, where
\begin{equation}\label{eq:phys_sol.2}
\begin{split}
& \overline{D}_t := \inf\{y>0:\,y - F_t(y)>0\} \\
& \quad\,:=\inf\bigg\{y>0:\,y + C\log\bigg(1-\frac{\PP(\overline{\tau}\ge t,\,\overline{Y}_{t-}\in(0,y))}
{\PP(\overline{\tau}\ge t)}\bigg)>0\bigg\}<\infty.
\end{split}
\end{equation}
\end{definition}

\begin{rmk}
It is easy to see that, for any $t\in[0,\overline{\tau}^0)\cap[0,T]$, it holds $\PP(\overline{\tau}> t)>0$, and, hence, all quantities in Definition \ref{def_phys_sol.new} are well-defined.
\end{rmk}


The path $\overline{Y}_t$, $t\in[0,\overline{\tau})\cap[0,T]$, for a physical solution $\overline{Y}$, should be thought of as the logarithmic asset value process of a \textit{typical} bank in a \textit{large} banking system, which is made precise by the following theorem. 

\begin{theorem}\label{main2}
Suppose that, for all $N\in\nn$, the initial values $Y_0^1,\,Y_0^2,\,\ldots,\,Y^N_0$ are i.i.d. according to a probability measure $\nu$ on $[0,\infty)$, with a bounded density $f_\nu$ vanishing in a neighborhood of $0$. Then, the sequence of random measures $\widetilde{\mu}^N$, $N\in\nn$ is tight with respect to the topology of weak convergence, and every limit point of this sequence belongs with probability one to the space of distributions of physical solutions $\overline{Y}$ for which $\overline{Y}_0\stackrel{d}{=}\nu$.
\end{theorem}

Note that Theorem \ref{main2}, in particular, proves the weak existence of a physical solution. The function $\Lambda$ in Definition \ref{def_phys_sol.new} represents the aggregate losses (on the logarithmic scale) of a typical bank in a large banking system resulting from the defaults of other banks. When default cascades lead to a jump in $\Lambda$, we speak of a \textit{systemic event}. At the \textit{random} time $\overline{\tau}$ the bank in consideration defaults; the \textit{deterministic} quantity $\overline{\tau}^0$ represents the time when the last bank defaults; $\overline{D}_t$ provides an upper bound on the maximum logarithmic value among the banks defaulting at time $t$; and $F_t(\overline{D}_t)$ is an upper bound on the total losses due to mutual exposure (on the logarithmic scale) incurred at time $t$. The term ``physical solution'' is borrowed from \cite{Delarue2}. 

\medskip

Our main interest is the time of the first \emph{systemic event}
\begin{equation}
t_{sys} := \inf\{t\in[0,\overline{\tau}^0)\cap[0,T]\,:\,\Lambda_t\neq\Lambda_{t-}\},
\end{equation}
when a non-negligible fraction of banks defaults in a short period of time. The time $t_{sys}$ can be viewed as the time of the first \emph{phase transition}, with the banking system passing abruptly from the well-behaved regime to the systemic crisis regime.

\medskip

Assuming that $\overline{Y}_0$ admits a density, we read off from Definition \ref{def_phys_sol.new} of a physical solution that 
\begin{equation}\label{tsysLBD}
t_{sys}\geq \inf\{t\in[0,\overline{\tau}^0)\cap[0,T]\,:\,r^*_t \geq1/C\},
\end{equation}
where 
\begin{equation}\label{r_star_t}
r^*_t:= \lim_{\eta\downarrow 0}\,\sup_{s\in(0,t]}\,\mathrm{ess\,sup}_{y\in(0,\eta)}\,\frac{p(s,y)}{\PP(\overline{\tau}\ge s)}
\end{equation}
and $p(s,\cdot)$ is the density of the distribution of $\overline{Y}_{s-}\,\mathbf{1}_{\{\overline{\tau}\ge s\}}$ restricted to $(0,\infty)$ (see Lemma \ref{le:phys.sol.1} below for the existence of $p(s,\cdot)$). 
The next theorem gives the corresponding upper bound on $t_{sys}$ in terms of the normalized density $p(t,y)/\PP(\overline{\tau}\ge t)$ near $y=0$:
\begin{equation}\label{tsysUBD}
t_{sys}\leq \inf\{t\in[0,\overline{\tau}^0)\cap[0,T]\,:\,r^{**}_t > c^*/C\},
\end{equation}
where 
\begin{equation}
r^{**}_t:= \lim_{\eta\downarrow 0}\,\mathrm{ess\,inf}_{y\in(0,\eta)}\,\frac{p(t,y)}{\PP(\overline{\tau}\ge t)},
\end{equation}
and $c^*$ is a constant depending only on $\sigma$.

\begin{theorem}\label{main2.2}
There exists a constant $c^*=c^*(\sigma)<\infty$ such that, for any c\`adl\`ag process $\overline{Y}$ satisfying (\ref{eq:phys_sol.1}), with the associated $\Lambda$, $\overline{\tau}^0$ and $r^{**}$, and with $\overline{Y}_0$ admitting a density, for any $t\in[0,\overline{\tau}^0)$, we have $\Lambda_{t}\neq\Lambda_{t-}$ whenever $r^{**}_t>c^*/C$. In particular, the upper bound \eqref{tsysUBD} holds.
\end{theorem}


The results of \cite{Delarue2} show that, for the model considered therein, one can take $c^*=1$, which makes it a natural to conjecture the same in the present case. However, the methods used in \cite{Delarue2} to establish this result do not seem to be applicable in the present case, and the methods used herein do not allow us to show that $c^*=1$.

\medskip

The inequalities \eqref{tsysLBD}, \eqref{tsysUBD} show that the normalized density $p(t,y)/\PP(\overline{\tau}\ge t)$ near $y=0$ can be used to \textit{measure the proximity} to the time $t_{sys}$ of the first systemic event. Simply put, a systemic event occurs when the normalized density at zero reaches the level $c^*/C$. This observation yields a natural connection between a systemic event and the two relevant \emph{observable} quantities: the fraction of banks at immediate risk and the level of \emph{non-core exposure}.
In addition, it motivates the analysis of physical solutions in relation to the values of the normalized density at zero. 
Our next result is concerned with the \textit{regularity} of a physical solution in relation to the normalized density at zero.
Lemma 5.2 in \cite{Delarue1} establishes the $1/2$-H\"older continuity of the cumulative loss process at any time at which it does not jump. The subsequent results in \cite{Delarue1} show that the normalized density at zero vanishes and the cumulative loss process becomes continuously differentiable at all times, if the strength of the interaction $C$ is sufficiently small (the model analyzed in \cite{Delarue1} is not exactly the same as the present one, but the arguments used therein can be adapted to the present case). The following theorem fills the gap between the two results: it shows that the cumulative loss process possesses higher H\"older regularity (even though it may not be continuously differentiable) if the normalized density at zero vanishes, without the assumption that $C$ is sufficiently small.

\begin{theorem}\label{thm:Holder.loss}
Let $\overline{Y}$ be a c\`adl\`ag process satisfying (\ref{eq:phys_sol.1}), with the associated $\Lambda$ and $\overline{\tau}^0$. Suppose that $\overline{Y}_0$ has a bounded density vanishing in a neighborhood of $0$. Consider any $t_0\in(0,\overline{\tau}^0)$ for which $r^{*}_{t_0}=0$. Then, for any $t_0'\in[0,t_0)$ there exist $\widetilde{C}<\infty$ and $\gamma\in(0,1]$ such that
\begin{equation}\label{eq.PhysSolReg.L.Holder}
|\Lambda_{t}-\Lambda_{s}| \leq \widetilde{C}\,|t-s|^{(1+\gamma)/2},\;\;s,t\in[0,t_0'].
\end{equation}
\end{theorem}

Next, we turn to the uniqueness of a physical solution. Note that establishing uniqueness is not only interesting in its own right, but it would also strengthen the convergence result significantly. Indeed, once the uniqueness is established, Theorem \ref{main2} would imply that $\widetilde{\mu}^N$, $N\in\mathbb{N}$ converge to a deterministic limit, which is the law of the unique physical solution. To date, the uniqueness of a general physical solution given by Definition \ref{def_phys_sol.new}, or its analogue in \cite{Delarue1}, \cite{Delarue2}, remains an open problem. Nevertheless, the uniqueness can be established in a class of sufficiently regular solutions, which can be described via an associated Cauchy-Dirichlet system. Such a uniqueness result is established in \cite{Delarue1}, under the additional assumption that $C$ is sufficiently small, which ensures that the cumulative loss process is continuously differentiable and, in particular, rules out the possibility of a jump. A local uniqueness result is also established in \cite{Delarue1}, and it does not require $C$ to be sufficiently small. Nevertheless, the latter result only holds on a time interval on which the cumulative loss process is continuously differentiable.
Herein, we do not make an assumption that $C$ is small, as we would like to analyze systems in which the cumulative loss process can jump. In addition, we establish uniqueness on a time interval on which the loss process neither jumps nor possesses a continuous derivative.
More specifically, we establish uniqueness up to the time
\begin{equation}\label{eq.treg.def}
t_{reg}=\big(\sup\,\big\{t\in(0,\overline{\tau}^0):\,\|\lambda\|_{L^2([0,t])}<\infty\big\}\big)\wedge T,
\end{equation}
where $\lambda$ is the weak derivative of $\Lambda$, and we use the conventions: $\sup \emptyset = 0$ and $\|\lambda\|_{L^2([0,t])}=\infty$ if $\Lambda$ is not absolutely continuous on $[0,t]$.
The following theorem proves the uniqueness of the \textit{stopped} physical solution $\overline{Y}_{t\wedge\overline{\tau}}$, $t\in[0,t_{reg})$, in the class of solutions with $t_{reg}>0$ and such that $\|\lambda\|_{L^2([0,\cdot])}$ ``does not jump to infinity". Moreover, it provides a precise connection between the cumulative loss process and the normalized density on $[0,t_{reg})$. 
\begin{theorem}\label{main1}
Let $\nu$ be a probability measure on $[0,\infty)$ admitting a density $f_\nu$ in the Sobolev space $W^1_2([0,\infty))$ with $f_\nu(0)=0$. Then,
\begin{enumerate}[(a)]
\item there exists a physical solution $\overline{Y}$, such that $\overline{Y}_0\stackrel{d}{=}\nu$ and the associated $t_{reg}$ and $\lambda$ satisfy
$$
t_{reg}>0,\quad \lim_{t\,\uparrow\, t_{reg}} \|\lambda\|_{L^2([0,t])} = \infty;
$$
\item the value of $t_{reg}>0$ is the same for all physical solutions satisfying the conditions of part (a), and the corresponding stopped physical solutions $\overline{Y}_{t\wedge\overline{\tau}}$, $t\in[0,t_{reg})$ are indistinguishable;
\item $p(\cdot,\cdot)$ is continuous on $[0,t_{reg})\times[0,\infty)$, with $p(\cdot,0)\equiv0$; moreover, the weak derivative $\partial_yp$ satisfies $(\partial_yp)(\cdot,0)\in L^2_{loc}([0,t_{reg}))$ and 
\begin{equation}\label{expl_crit}
\lambda_t=-C\,\frac{\sigma^2}{2}\,\frac{(\partial_yp)(t,0)}{\int_0^{\infty} p(t,y)dy}\;\;\text{for\;almost\;every\;}\;t\in[0,t_{reg}).
\end{equation} 
\end{enumerate}
\end{theorem}



Parts (a) and (b) of Theorem \ref{main1} show that the logarithmic asset value of a typical bank in a large banking system behaves according to the unique stopped physical solution until the time $t_{reg}>0$, given by \eqref{eq.treg.def}.
Theorem \ref{main1}(c) expresses the value of $\lambda_t$ through the slope of the normalized logarithmic asset value profile of banks that are close to failure at time $t$. 
\medskip

The rest of the paper is structured as follows. In Section \ref{sec:cdp} we analyze the Cauchy-Dirichlet problem associated with a stopped physical solution before the explosion of $\lambda$ in the $L^2$ norm, which is used in the proof of Theorem \ref{main1}. Section \ref{sec:rps} studies the fixed-point problem satisfied by $\lambda$ until it explodes in the $L^2$ norm. We use Sobolev norm estimates for solutions to linear parabolic PDEs in \cite{Lady} (see \cite[Chapter III, Section 6]{LSU}) and parabolic Sobolev inequalities (see e.g. \cite[Chapter II, Lemmas 3.3, 3.4]{LSU}) to show that the Banach fixed-point theorem is applicable to a suitable ``truncated'' fixed-point problem. This yields the existence and uniqueness of the solution to the original fixed-point problem. The latter is used to construct the unique stopped physical solution until the explosion of $\lambda$ in the $L^2$ norm, proving Theorem \ref{main1}. Section \ref{se:regularity} establishes a priori regularity properties of physical solutions and connects the behavior of the normalized density $p(t,y)/\PP(\overline{\tau}\ge t)$ near $y=0$ with the H\"older continuity of $\Lambda$, proving Theorem \ref{thm:Holder.loss}. Section \ref{se:sufficientCond.Jump} provides the proof of Theorem \ref{main2.2}. Section \ref{se:convergence} is devoted to the proof of Theorem \ref{main2}, which is an adaptation of the arguments used in \cite{Delarue2}.

\section{Cauchy-Dirichlet problem} \label{sec:cdp}

For $\nu$ as in Theorem \ref{main1}, $T_1\in(0,\infty)$, and $\lambda\in L^2([0,T_1])$ consider the Cauchy-Dirichlet problem 
\begin{equation} \label{eq:cdp}
\partial_t p = -(\alpha+\lambda_t)\,\partial_y p + \frac{\sigma^2}{2}\,\partial_y^2 p,\;\; p(0,\cdot)=f_\nu,\;\; p(\cdot,0)=0.  
\end{equation}
The next two lemmas investigate its solution $p$. 

\begin{lemma}\label{lemma:exp_bd}
Let $\nu$ be as in Theorem \ref{main1}. Then, for $T_1\in(0,\infty)$ and $\lambda\in L^2([0,T_1])$, there exists a unique generalized solution $p$ of \eqref{eq:cdp} in the space $W^{1,2}_2([0,T_1]\times[0,\infty))$. Moreover, $p$ is non-negative and satisfies the integrability estimates
\begin{equation}\label{eq:exp_est}
\begin{split}
& e^{-\int_0^t (\frac{\alpha+\lambda_s}{\sigma})^2\,\mathrm{d}s} \int_0^\infty 
\bigg(2\Phi\bigg(\frac{y}{\sigma\sqrt{t}}\bigg)-1\bigg)^2\,f_\nu(y)\,\mathrm{d}y \\
& \le \int_0^\infty p(t,y)\,\mathrm{d}y 
\le e^{\frac{1}{2} \int_0^t (\frac{\alpha+\lambda_s}{\sigma})^2\,\mathrm{d}s} \int_0^\infty 
\bigg(2\Phi\bigg(\frac{y}{\sigma\sqrt{t}}\bigg)-1\bigg)^{1/2}
\,f_\nu(y)\,\mathrm{d}y
\end{split}
\end{equation}
for all $t\in[0,T_1]$, where $\Phi$ is the standard Gaussian cumulative distribution function.
\end{lemma}

\begin{lemma}\label{lemma:exp_test}
Let $\nu$ be as in Theorem \ref{main1}. Then, for $T_1\in(0,\infty)$ and $\lambda\in L^2([0,T_1])$, the unique generalized solution $p\in W^{1,2}_2([0,T_1]\times[0,\infty))$ of \eqref{eq:cdp} fulfills
\begin{equation}\label{eq:exp_test}
\int_0^\infty p(t,y)\,\mathrm{d}y
=\int_0^\infty f_\nu(y)\,\mathrm{d}y
-\frac{\sigma^2}{2}\,\int_0^t(\partial_y p)(s,0)\,\mathrm{d}s,\;\;t\in[0,T_1].
\end{equation}
\end{lemma}

\smallskip

\noindent\textbf{Proof of Lemma \ref{lemma:exp_bd}. Step 1.} The existence and uniqueness of the generalized solution $p\in W^{1,2}_2([0,T_1]\times[0,\infty))$ of \eqref{eq:cdp} follow from the results of \cite[Chapter III, Section 6]{LSU} (see \cite[Chapter III, Remark 6.3]{LSU} and note that $-(\alpha+\lambda)$ fulfills the condition (6.26) there). We also refer to the original reference \cite{Lady}. 

\medskip

Now, consider the process $\overline{Z}^\lambda_t$, $t\in[0,T_1]$ defined by
\begin{equation}\label{eq:Ybar_def}
\begin{split}
& \overline{Z}^\lambda_0\stackrel{d}{=}\nu, \quad
\overline{Z}^\lambda_t=\overline{Z}^\lambda_0+\alpha\,t
+\int_0^t \lambda_s\,\mathrm{d}s + \sigma\,\overline{B}_t,\;\; t\in[0,T_1\wedge\overline{\tau}],\\
& \overline{\tau}=\inf\{t\in[0,T_1]:\,\overline{Z}^\lambda_t=0\},\quad
\overline{Z}^\lambda_t=0,\;\; \overline{\tau}\le t\le T_1,
\end{split}
\end{equation}
where $\overline{B}$ is a standard Brownian motion independent of $\overline{Z}^\lambda_0\stackrel{d}{=}\nu$. The Radon-Nikodym and the Girsanov theorems show that the law of $\overline{Z}^\lambda_t$ has a density with respect to that of $(\overline{Z}_0^\lambda+\sigma\,\overline{B})_{t\wedge\overline{\tau}}$ for all $t\in[0,T_1]$. In particular, the restriction of the law of $\overline{Z}^\lambda_t$ to $(0,\infty)$ possesses a density $\widetilde{p}(t,\cdot)$ with respect to the Lebesgue measure for all $t\in[0,T_1]$. We claim next that the $W^{1,2}_2([0,T_1]\times[0,\infty))$-solution $p$ of \eqref{eq:cdp} equals $\widetilde{p}$. 

\medskip

\noindent\textbf{Step 2.} We fix a $t\in[0,T_1]$, pick a function $h\in W^1_2([0,\infty))$ with $h(0)=0$, and consider the auxiliary problem  
\begin{equation}\label{eq:zeta_problem}
\partial_s\zeta+(\alpha+\lambda_s)\,\partial_y\zeta+\frac{\sigma^2}{2}\,\partial_y^2\zeta=0, \;\;\zeta(t,\cdot)=h,\;\;\zeta(\cdot,0)=0,\;\;\zeta\in W^{1,2}_2([0,t]\times[0,\infty)).
\end{equation} 
As with the problem \eqref{eq:cdp} there exists a unique generalized solution $\zeta$ of \eqref{eq:zeta_problem}. Moreover, for any fixed $K\in(0,\infty)$ and with 
\begin{equation}
\overline{\tau}_K:=
\inf\{s\in[0,T_1]:\,\overline{Z}^\lambda_0+\sigma\,\overline{B}_s=K\}
\end{equation}
the PDE in \eqref{eq:zeta_problem} and the It\^o formula in \cite[Section 2.10, Theorem 1]{Krylov2} yield
\begin{equation}\label{zeta_Ito}
\begin{split}
\zeta(t\wedge\overline{\tau}\wedge\overline{\tau}_K,(\overline{Z}^\lambda_0+\sigma\,\overline{B})_{t\wedge\overline{\tau}\wedge\overline{\tau}_K})
=\zeta(0,\overline{Z}^\lambda_0)+\int_0^{t\wedge\overline{\tau}\wedge\overline{\tau}_K} (\partial_y\zeta)(s,\overline{Z}^\lambda_0+\sigma\,\overline{B}_s)\,\sigma\,\mathrm{d}\overline{B}_s \\
-\int_0^{t\wedge\overline{\tau}\wedge\overline{\tau}_K}
(\partial_y\zeta)(s,\overline{Z}^\lambda_0+\sigma\,\overline{B}_s)\,(\alpha+\lambda_s)\,\mathrm{d}s
\end{split}
\end{equation}
(note that $\zeta$ is a continuous bounded function and $(\partial_y\zeta)\in L^4([0,t]\times[0,\infty))$ by the parabolic Sobolev inequality in the form of \cite[Chapter II, Lemma 3.3]{LSU}). In view of the Girsanov Theorem, \eqref{zeta_Ito} implies
\begin{equation}\label{zeta_Ito2}
\zeta(t\wedge\overline{\tau}\wedge\overline{\tau}_K,\overline{Z}^\lambda_{t\wedge\overline{\tau}\wedge\overline{\tau}_K})
=\zeta(0,\overline{Z}^\lambda_0)+\int_0^{t\wedge\overline{\tau}\wedge\overline{\tau}_K} (\partial_y\zeta)(s,\overline{Z}^\lambda_s)\,\sigma\,\mathrm{d}\overline{B}_s.
\end{equation}

\smallskip

Next, we combine the Girsanov theorem with H\"older's and Jensen's inequalities to obtain the chain of estimates
\begin{equation*}
\begin{split}
& \E\bigg[\int_0^{t\wedge\overline{\tau}\wedge\overline{\tau}_K} \big((\partial_y\zeta)(s,\overline{Z}^\lambda_s)\big)^2\,\mathrm{d}s\bigg] \\
& =\E\bigg[e^{-\int_0^t \frac{\alpha+\lambda_s}{\sigma}\,\mathrm{d}\overline{B}_s-\frac{1}{2}\int_0^t
(\frac{\alpha+\lambda_s}{\sigma})^2\,\mathrm{d}s}\int_0^{t\wedge\overline{\tau}\wedge\overline{\tau}_K} \big((\partial_y\zeta)(s,\overline{Z}^\lambda_0+\sigma\,\overline{B}_s)\big)^2\,\mathrm{d}s\bigg] \\
& \le\E\bigg[e^{-3\int_0^t \frac{\alpha+\lambda_s}{\sigma}\,\mathrm{d}\overline{B}_s-\frac{3}{2}\int_0^t
(\frac{\alpha+\lambda_s}{\sigma})^2\,\mathrm{d}s}\bigg]^{1/3}\,
\E\bigg[\bigg(\int_0^{t\wedge\overline{\tau}\wedge\overline{\tau}_K} 
\big((\partial_y\zeta)(s,\overline{Z}^\lambda_0+\sigma\,\overline{B}_s)\big)^2\,\mathrm{d}s\bigg)^{3/2}\bigg]^{2/3} \\
& \le e^{\int_0^t (\frac{\alpha+\lambda_s}{\sigma})^2\,\mathrm{d}s}\,t^{1/3}\,\E\bigg[\int_0^{t\wedge\overline{\tau}\wedge\overline{\tau}_K} \big|(\partial_y\zeta)(s,\overline{Z}^\lambda_0+\sigma\,\overline{B}_s)\big|^3\,\mathrm{d}s
\bigg]^{2/3}.
\end{split}
\end{equation*}
The latter expression is finite thanks to \cite[Section 2.2, Theorem 4]{Krylov2} and $(\partial_y\zeta)\in L^6([0,t]\times[0,\infty))$ (a consequence of the parabolic Sobolev inequality in the form of \cite[Chapter II, Lemma 3.3]{LSU}). Consequently, taking the expectation on both sides of \eqref{zeta_Ito2} and passing to the limit $K\to\infty$ we get 
\begin{equation}
\E[\zeta(t\wedge\overline{\tau},\overline{Z}^\lambda_{t\wedge\overline{\tau}})]=\E[\zeta(0,\overline{Z}^\lambda_0)],
\end{equation}
which can be rewritten as 
\begin{equation}\label{tilde_p_mart}
\int_0^\infty h(y)\,\widetilde{p}(t,y)\,\mathrm{d}y
=\int_0^\infty \zeta(0,y)\,f_\nu(y)\,\mathrm{d}y.
\end{equation}

\smallskip

On the other hand, $p\in W^{1,2}_2([0,T_1]\times[0,\infty))$ implies that the norms $\|p(s,\cdot)\|_{L^2([0,\infty))}$, $s\in[0,t]$ are uniformly bounded due to the continuity of the evaluation map (see e.g. \cite[Chapter II, Lemma 3.4]{LSU}). This and a density argument invoking the continuity of the evaluation map one more time show that the weak formulation of the problem \eqref{eq:cdp} applies to test functions in $W^{1,2}_2([0,t]\times[0,\infty))$. For the solution $\zeta$ of \eqref{eq:zeta_problem} it gives
\begin{equation}
\int_0^\infty h(y)\,p(t,y)\,\mathrm{d}y
=\int_0^\infty \zeta(0,y)\,f_\nu(y)\,\mathrm{d}y,
\end{equation}
which together with \eqref{tilde_p_mart} and the arbitrariness of $h$, $t$ implies $\widetilde{p}=p$ on $[0,T_1]\times[0,\infty)$.

\medskip

\noindent\textbf{Step 3.} The non-negativity of $p$ is now an immediate consequence of the non-negativity of $\widetilde{p}$. In addition, $\int_0^\infty p(t,y)\,\mathrm{d}y$ can be rewritten as 
\begin{equation*}
\E\big[\mathbf{1}_{\{\overline{Z}^\lambda_t>0\}}\big]
=\int_0^\infty \E\Big[e^{-\int_0^t \frac{\alpha+\lambda_s}{\sigma}\,\mathrm{d}\overline{B}_s-\frac{1}{2}\int_0^t
(\frac{\alpha+\lambda_s}{\sigma})^2\,\mathrm{d}s}\,\mathbf{1}_{\{y+\sigma\overline{B}_s>0,\,0\le s\le t\}}\Big]
\,f_\nu(y)\,\mathrm{d}y.
\end{equation*}
At this point, the estimates of \eqref{eq:exp_est} follow from the Cauchy-Schwarz inequality in the forms 
\begin{equation*}
\begin{split}
& \E\Big[e^{\int_0^t \frac{\alpha+\lambda_s}{\sigma}\,\mathrm{d}\overline{B}_s
+\frac{1}{2}\int_0^t (\frac{\alpha+\lambda_s}{\sigma})^2\,\mathrm{d}s}\Big]^{-1}\,
\pp\Big(\overline{B}_s>-\frac{y}{\sigma},\;0\le s\le t\Big)^2 \\
& \le\E\Big[e^{-\int_0^t \frac{\alpha+\lambda_s}{\sigma}\,\mathrm{d}\overline{B}_s-\frac{1}{2}\int_0^t
(\frac{\alpha+\lambda_s}{\sigma})^2\,\mathrm{d}s}\,\mathbf{1}_{\{y+\sigma\overline{B}_s>0,\;0\le s\le t\}}\Big] \\
& \le \E\Big[e^{-2\int_0^t \frac{\alpha+\lambda_s}{\sigma}\,\mathrm{d}\overline{B}_s
-\int_0^t (\frac{\alpha+\lambda_s}{\sigma})^2\,\mathrm{d}s}\Big]^{1/2}\,\pp\Big(\overline{B}_s>-\frac{y}{\sigma},\;0\le s\le t\Big)^{1/2}
\end{split}
\end{equation*}
and the reflection principle for Brownian motion. \ep

\medskip

We proceed to the proof of Lemma \ref{lemma:exp_test}.  

\medskip

\noindent\textbf{Proof of Lemma \ref{lemma:exp_test}.} We pick a sequence $h_n$, $n\in\nn$ of infinitely differentiable functions on $[0,\infty)$ such that
\begin{enumerate}[(i)]
\item $h_n(y)=1$ if $n^{-1}\le y\le n$ and $h_n(y)=0$ if $y\le(n+1)^{-1}$ or $y\ge n+1$,
\item $h_n$ is non-decreasing on $[(n+1)^{-1},n^{-1}]$ and non-increasing on $[n,n+1]$,
\item $\sup_{n\in\nn}\,\sup_{[n,n+1]} |h_n'|<\infty$ and $\sup_{n\in\nn}\,\sup_{[n,n+1]} |h_n''|<\infty$.   
\end{enumerate}
The weak formulation of \eqref{eq:cdp} for each such function reads
\begin{equation}\label{weak_form_approx}
\begin{split}
& \int_0^\infty h_n(y)\,p(t,y)\,\mathrm{d}y
-\int_0^\infty h_n(y)\,f_\nu(y)\,\mathrm{d}y \\
& =\int_0^t (\alpha+\lambda_s) \int_0^\infty h_n'(y)\,p(s,y)\,\mathrm{d}y\,\mathrm{d}s +\frac{\sigma^2}{2} \int_0^t \int_0^\infty h_n''(y)\,p(s,y)\,\mathrm{d}y\,\mathrm{d}s,\;\;t\in[0,T_1].
\end{split}
\end{equation}

\smallskip

The monotone convergence theorem implies that the first line in \eqref{weak_form_approx} tends to
\begin{equation*}
\int_0^\infty p(t,y)\,\mathrm{d}y - \int_0^\infty f_\nu(y)\,\mathrm{d}y
\end{equation*}
in the limit $n\to\infty$. Moreover, the first summand on the second line in \eqref{weak_form_approx} can be rewritten as 
\begin{equation}\label{weak_form_approx1}
\int_0^t (\alpha+\lambda_s) \int_{(n+1)^{-1}}^{n^{-1}} h_n'(y)\,p(s,y)\,\mathrm{d}y\,\mathrm{d}s
+\int_0^t (\alpha+\lambda_s) \int_n^{n+1} h_n'(y)\,p(s,y)\,\mathrm{d}y\,\mathrm{d}s.
\end{equation}
Combining $p(s,0)=0$, $s\in[0,T]$, the uniform continuity of $p$ on $[0,T_1]\times[0,1]$ (due to the parabolic Sobolev inequality in \cite[Chapter II, Lemma 3.3]{LSU}), and property (ii) above we see that the first summand in \eqref{weak_form_approx1} converges to $0$ as $n\to\infty$. The same is true for the second summand in \eqref{weak_form_approx1} thanks to property (iii) above, the upper bound of Lemma \ref{lemma:exp_bd}, and the dominated convergence theorem. 

\medskip

The second summand on the second line in \eqref{weak_form_approx} can be recast as 
\begin{equation}\label{weak_form_approx2}
\frac{\sigma^2}{2} \int_0^t \int_n^{n+1} h_n''(y)\,p(s,y)\,\mathrm{d}y\,\mathrm{d}s \\
+\frac{\sigma^2}{2} \int_0^t \int_{(n+1)^{-1}}^{n^{-1}} h_n''(y)\,p(s,y)\,\mathrm{d}y\,\mathrm{d}s.
\end{equation}
As $n\to\infty$, the first summand in \eqref{weak_form_approx2} converges to $0$ by the same argument as used to analyze the second summand in \eqref{weak_form_approx1}. Next, we employ integration by parts to transform the second summand in \eqref{weak_form_approx2} to 
\begin{equation}\label{weak_form_approx3}
-\frac{\sigma^2}{2} \int_0^t \int_{(n+1)^{-1}}^{n^{-1}} h_n'(y)\,(\partial_y p)(s,y)\,\mathrm{d}y\,\mathrm{d}s
\end{equation}
(recall that $p(s,\cdot)\in W^1_2([0,\infty))$, $s\in[0,T_1]$ thanks to the well-definedness of the first evaluation map in \cite[Chapter II, Lemma 3.4]{LSU}). The quantity in \eqref{weak_form_approx3} converges to $-\frac{\sigma^2}{2}\,\int_0^t (\partial_y p)(s,0)\,\mathrm{d}s$ as $n\to\infty$, since
\begin{equation*}
\begin{split}
& \limsup_{n\to\infty} \bigg|\int_0^t \int_{(n+1)^{-1}}^{n^{-1}} h_n'(y)\,(\partial_y p)(s,y)\,\mathrm{d}y\,\mathrm{d}s
-\int_0^t (\partial_y p)(s,0)\,\mathrm{d}s\bigg| \\
& =\limsup_{n\to\infty} \bigg|\int_0^t \int_{(n+1)^{-1}}^{n^{-1}}
h_n'(y)\,\int_0^y (\partial_y^2 p)(s,z)\,\mathrm{d}z\,\mathrm{d}y\,\mathrm{d}s\bigg| \\
& \le\limsup_{n\to\infty} \int_0^t \int_0^{n^{-1}} |\partial_y^2 p|(s,z)\,\mathrm{d}z\,\mathrm{d}s=0,
\end{split}
\end{equation*}
where we have relied on properties (i), (ii) above and $p\in W^{1,2}_2([0,T_1]\times[0,\infty))$. All in all, we end up with \eqref{eq:exp_test} when we take the $n\to\infty$ limit in \eqref{weak_form_approx}. \ep

\section{Regular interval of the physical solution} \label{sec:rps}

The next proposition is the key to the proof of Theorem \ref{main1} and establishes the existence and uniqueness of the solution to the fixed-point problem associated with the function $\lambda$ in \eqref{eq.treg.def}. 

\begin{proposition}\label{prop:ex&uniq}
Let $\nu$ be as in Theorem \ref{main1}. Then, \smallskip
\begin{enumerate}[(a)]
\item there exist a time $t_{reg}\in(0,T]$ and a function $\lambda\in L^2_{loc}([0,t_{reg}))$ such that for all $T_1\in(0,t_{reg})$ the unique generalized solution of 
\begin{equation}\label{FPprob1}
\partial_t p = -(\alpha+\lambda_t)\,\partial_y p + \frac{\sigma^2}{2}\,\partial_y^2 p,\;\; p(0,\cdot)=f_\nu,\;\; p(\cdot,0)=0
\end{equation}
in $W^{1,2}_2([0,T_1]\times[0,\infty))$ satisfies
\begin{equation}\label{FPprob2}
-C\,\frac{\sigma^2}{2}\,\frac{(\partial_y p)(t,0)}{\int_0^\infty p(t,y)\,\mathrm{d}y}=\lambda_t\;\;\text{for\;almost\;every}\;\;t\in[0,T_1]
\end{equation}
and $\lim_{T_1\uparrow t_{reg}} \|\lambda\|_{L^2([0,T_1])}=\infty$ if $t_{reg}<T$; \medskip
\item for any $(t_{reg},\lambda)$, $(\widetilde{t}_{reg},\widetilde{\lambda})$ such that 
\begin{enumerate}[(i)]
\item $t_{reg},\widetilde{t}_{reg}\in(0,T]$, 
\item $\lambda\in L^2_{loc}([0,t_{reg}))$, $\widetilde{\lambda}\in L^2_{loc}([0,\widetilde{t}_{reg}))$ satisfy the fixed-point problem \eqref{FPprob1}, \eqref{FPprob2} for all $T_1\in(0,t_{reg})$, $T_1\in(0,\widetilde{t}_{reg})$, respectively, 
\item $\lim_{T_1\uparrow t_{reg}} \|\lambda\|_{L^2([0,T_1])}=\infty$ if $t_{reg}<T$, $\lim_{T_1\uparrow \widetilde{t}_{reg}} \|\widetilde{\lambda}\|_{L^2([0,T_1])}=\infty$ if $\widetilde{t}_{reg}<T$ 
\end{enumerate}
it holds $t_{reg}=\widetilde{t}_{reg}$ and $\lambda=\widetilde{\lambda}$ almost everywhere.
\end{enumerate}
\end{proposition}

\medskip

\noindent\textbf{Proof of Proposition \ref{prop:ex&uniq}. Step 1.} Our first aim is to show for all $M\in(0,\infty)$ and all small enough $T_1=T_1(M)\in(0,T)$ the existence and uniqueness in $L^2([0,T_1])$ for the ``truncated'' fixed-point problem
\begin{equation}
\partial_t p=-(\alpha+\lambda^{M,T_1}_t)\,\partial_y p+\frac{\sigma^2}{2}\,\partial^2_y p,\;\; p(0,\cdot)=f_\nu,\;\; p(\cdot,0)=0,\;\; p\in W^{1,2}_2([0,T_1]\times[0,\infty)), 
\label{Krylov_tr1}
\end{equation}
\begin{equation}
-C\,\frac{\sigma^2}{2}\,\frac{(\partial_y p)(t,0)}{\int_0^\infty p(t,y)\,\mathrm{d}y}=\lambda_t\;\;\text{for\;almost\;every}\;\;t\in[0,T_1],\qquad\qquad\qquad\qquad\qquad\qquad\qquad
\label{Krylov_tr2}
\end{equation}
where 
\begin{equation}
\lambda^{M,T_1}:=\lambda\,\mathbf{1}_{\{\|\lambda\|_{L^2([0,T_1])}\le M\}}+\lambda\,\frac{M}{\|\lambda\|_{L^2([0,T_1])}}\,\mathbf{1}_{\{\|\lambda\|_{L^2([0,T_1])}>M\}}. 
\end{equation}
To this end, it suffices to verify that the mapping taking $L^2([0,T_1])$ functions $\lambda$ to the left-hand side of \eqref{Krylov_tr2} (with $p$ being the unique generalized solution of \eqref{Krylov_tr1}) is a contraction on $L^2([0,T_1])$, since then the Banach fixed-point theorem can be applied. We observe that the described mapping is well-defined with its range contained in $L^2([0,T_1])$ by the assumptions $f_\nu\in W^1_2([0,\infty))$, $f_\nu(0)=0$, the existence and uniqueness result of \cite[Chapter III, Remark 6.3]{LSU}, the well-definedness of the second evaluation map in \cite[Chapter II, Lemma 3.4]{LSU}, and the lower bound in \eqref{eq:exp_est}. The following two steps are devoted to the proof of the contraction property.

\medskip

\noindent\textbf{Step 2.} Given two $L^2([0,T_1])$ functions $\lambda$, $\widetilde{\lambda}$, let $p$, $\widetilde{p}$ be the corresponding solutions of \eqref{Krylov_tr1} and note that $\Delta:=p-\widetilde{p}\in W^{1,2}_2([0,T_1]\times[0,\infty))$ satisfies 
\begin{equation}\label{Delta_PDE}
\begin{split}
\partial_t\Delta=-(\alpha+\widetilde{\lambda}^{M,T_1}_t)\,\partial_y\Delta+\frac{\sigma^2}{2}\,\partial_y^2\Delta+(\widetilde{\lambda}^{M,T_1}_t-\lambda^{M,T_1}_t)\,\partial_yp,\;\; 
\Delta(0,\cdot)=0,\;\;\Delta(\cdot,0)=0.
\end{split}
\end{equation}
The source term in \eqref{Delta_PDE} admits the norm bound
\begin{equation}\label{source1_bound}
\begin{split}
& \big\|(\widetilde{\lambda}^{M,T_1}_t-\lambda^{M,T_1}_t)\,\partial_yp\big\|_{L^2([0,T_1]\times[0,\infty))} \\
& \le \big\|\widetilde{\lambda}^{M,T_1}_t-\lambda^{M,T_1}_t\big\|_{L^2([0,T_1])}\;\mathrm{ess\,sup}_{t\in[0,T_1]} \|(\partial_y p)(t,\cdot)\|_{L^2([0,\infty))} \\
& \le 2\,\|\widetilde{\lambda}_t-\lambda_t\big\|_{L^2([0,T_1])}\;\mathrm{ess\,sup}_{t\in[0,T_1]}\|(\partial_y p)(t,\cdot)\|_{L^2([0,\infty))}.
\end{split}
\end{equation}
Moreover, the boundedness of the first evaluation map in \cite[Chapter II, Lemma 3.4]{LSU} and the results of \cite[Chapter III, Section 6]{LSU} used for the solution $p$ of \eqref{Krylov_tr1} give the respective estimates 
\begin{equation}\label{ux_bound}
\mathrm{ess\,sup}_{t\in[0,T_1]} \|(\partial_y p)(t,\cdot)\|_{L^2([0,\infty))}\le C_1\,\|p\|_{W^{1,2}_2([0,T_1]\times[0,\infty))}\le C_2, 
\end{equation}
with constants $C_1=C_1(T)<\infty$ and $C_2=C_2(\alpha,M,\sigma,\|f_\nu\|_{W^1_2([0,\infty))},T)<\infty$. In view of \eqref{source1_bound}, \eqref{ux_bound}, we can now apply the boundedness of the first evaluation map in \cite[Chapter II, Lemma 3.4]{LSU} and the results of \cite[Chapter III, Section 6]{LSU} to the solution $\Delta$ of \eqref{Delta_PDE} to find 
\begin{equation}\label{Deltax_bound}
\mathrm{ess\,sup}_{t\in[0,T_1]} \|(\partial_y \Delta)(t,\cdot)\|_{L^2([0,\infty))}\le C_1\,\|\Delta\|_{W^{1,2}_2([0,T_1]\times[0,\infty))}\le C_3\,C_2\,\|\widetilde{\lambda}_t-\lambda_t\|_{L^2([0,T_1])},
\end{equation}
where the constant $C_3<\infty$ can be chosen in terms of $\alpha$, $M$, $\sigma$, and $T$ only.

\medskip

Next, we regard the PDE in \eqref{Delta_PDE} as a heat equation with the $L^2([0,T_1]\times[0,\infty))$ source 
\begin{equation}
g:=-(\alpha+\widetilde{\lambda}^{M,T_1}_t)\,\partial_y\Delta+(\widetilde{\lambda}^{M,T_1}_t-\lambda^{M,T_1}_t)\,\partial_yp.
\end{equation} 
In particular, we can write
\begin{eqnarray}
&&\quad\;\;\Delta(t,y)=\int_0^t \int_0^\infty g(s,z)\,\psi_\sigma(t-s,z,y)\,\mathrm{d}z\,\mathrm{d}s,\;\;(t,y)\in[0,T_1]\times[0,\infty),\\
&&\quad\;\;(\partial_y\Delta)(t,y)=\int_0^t \int_0^\infty g(s,z)\,(\partial_y\psi_\sigma)(t-s,z,y)\,\mathrm{d}z\,\mathrm{d}s,\;\;(t,y)\in[0,T_1]\times[0,\infty),
\end{eqnarray}
where 
\begin{equation}
\psi_\sigma(t-s,y,z):=(2\pi\sigma^2(t-s))^{-1/2}\bigg(\!\exp\bigg(\!-\frac{(y-z)^2}{2\sigma^2(t-s)}\bigg)
-\exp\bigg(\!-\frac{(y+z)^2}{2\sigma^2(t-s)}\bigg)\!\bigg)
\end{equation}
is the Dirichlet heat kernel on $[0,\infty)$ with the diffusion coefficient $\sigma$. It now follows from Fubini's theorem, the triangle inequality, Young's inequality, Cauchy-Schwarz inequality and \eqref{source1_bound}, \eqref{ux_bound}, \eqref{Deltax_bound} that 
\begin{equation}\label{short_time_heat_kernel}
\begin{split}
\|(\partial_y\Delta)(\cdot,0)\|_{L^2([0,T_1])}
&=\bigg\|\int_0^\infty \int_0^t g(s,z)\,(\partial_y\psi_\sigma)(t-s,z,0)\,\mathrm{d}s\,\mathrm{d}z\bigg\|_{L^2([0,T_1])} \\
&\le\int_0^\infty\bigg\|\int_0^t g(s,z)\,
(\partial_y\psi_\sigma)(t-s,z,0)\,\mathrm{d}s\bigg\|_{L^2([0,T_1])}\,\mathrm{d}z \\
&\le\int_0^\infty \|g(\cdot,z)\|_{L^2([0,T_1])}\,\|(\partial_y\psi_\sigma)(\cdot,z,0)\|_{L^1([0,T_1])}\,\mathrm{d}z \\
&\le\|g\|_{L^2([0,T_1]\times[0,\infty))}\,\bigg(\int_0^\infty \|(\partial_y\psi_\sigma)(\cdot,z,0)\|_{L^1([0,T_1])}^2\,\mathrm{d}z\bigg)^{1/2} \\
&\le C_4\,\|\widetilde{\lambda}_t-\lambda_t\|_{L^2([0,T_1])}\,
\bigg(\int_0^\infty \|(\partial_y\psi_\sigma)(\cdot,z,0)\|_{L^1([0,T_1])}^2\,\mathrm{d}z\bigg)^{1/2},
\end{split}
\end{equation}
with a constant $C_4=C_4(\alpha,M,\sigma,\|f_\nu\|_{W^1_2([0,\infty))},T)<\infty$.

\medskip

\noindent\textbf{Step 3.} Next, we subtract \eqref{eq:exp_test} for $\widetilde{p}$ from \eqref{eq:exp_test} for $p$, apply the triangle and the Cauchy-Schwarz inequalities, and use \eqref{short_time_heat_kernel} to find
\begin{equation}\label{eq:exp_Delta}
\begin{split}
\sup_{t\in[0,T_1]} \bigg|\int_0^\infty \Delta(t,y)\,\mathrm{d}y\bigg|
= \frac{\sigma^2}{2}\,\sup_{t\in[0,T_1]}\,\bigg|\int_0^t (\partial_y\Delta)(s,0)\,\mathrm{d}s\bigg| 
\le \frac{\sigma^2}{2}\,T_1^{1/2}\,\|(\partial_y\Delta)(\cdot,0)\|_{L^2([0,T_1])} \\
\le C_5\,T_1^{1/2}\,\|\widetilde{\lambda}_t-\lambda_t\|_{L^2([0,T_1])},
\end{split}
\end{equation}
where the constant $C_5<\infty$ depends on $\alpha$, $M$, $\sigma$, $f_\nu$, and $T$ only. 

\medskip

In addition, the triangle inequality and the lower bound in \eqref{eq:exp_est} imply
\begin{equation*}
\begin{split}
& \bigg\|\frac{(\partial_yp)(\cdot,0)}{\int_0^\infty p(\cdot,y)\,\mathrm{d}y}-\frac{(\partial_y\widetilde{p})(\cdot,0)}{\int_0^\infty \widetilde{p}(\cdot,y)\,\mathrm{d}y}\bigg\|_{L^2([0,T_1])} \\
& \le\bigg\|\frac{1}{\int_0^\infty p(\cdot,y)\,\mathrm{d}y}\,(\partial_y\Delta)(\cdot,0)\bigg\|_{L^2([0,T_1])} \qquad\\
&\quad +\bigg\|(\partial_y\widetilde{p})(\cdot,0)\,\bigg(\frac{1}{\int_0^\infty p(\cdot,y)\,\mathrm{d}y}-\frac{1}{\int_0^\infty \widetilde{p}(\cdot,y)\,\mathrm{d}y}\bigg)\bigg\|_{L^2([0,T_1])} \\
& \le C_6\,\bigg(\|(\partial_y\Delta)(\cdot,0)\|_{L^2([0,T_1])} 
+\bigg\|(\partial_y\widetilde{p})(\cdot,0)\,\int_0^\infty \Delta(\cdot,y)\,\mathrm{d}y\bigg\|_{L^2([0,T_1])}\bigg),
\end{split}
\end{equation*}
with a constant $C_6=C_6(\alpha,M,\sigma,f_\nu,T)<\infty$. In view of \eqref{short_time_heat_kernel}, \eqref{eq:exp_Delta}, and the boundedness of the second evaluation map in \cite[Chapter II, Lemma 3.4]{LSU}, the latter upper bound is at most
\begin{equation}\label{eq:contraction}
C_7\,\bigg(\bigg(\int_0^\infty \|(\partial_y\psi_\sigma)(\cdot,z,0)\|_{L^1([0,T_1])}^2\,\mathrm{d}z\bigg)^{1/2}+T_1^{1/2}\bigg)\,\|\widetilde{\lambda}_t-\lambda_t\|_{L^2([0,T_1])},
\end{equation}
where $C_7<\infty$ can be chosen in terms of $\alpha$, $M$, $\sigma$, $f_\nu$, and $T$ only. The desired contraction property for small enough $T_1=T_1(M)\in(0,T)$ readily follows. 

\medskip

\noindent\textbf{Step 4.} Now, we let 
\begin{equation}\label{max_exist_int}
t_{reg}:=\sup\{T_1\in(0,T):\text{the problem \eqref{FPprob1}, \eqref{FPprob2} has a solution }\lambda\!\in \!L^2([0,T_1])\}
\end{equation}
and claim that the supremum is taken over a non-empty set. Indeed, for fixed $M\in(0,\infty)$ and a small enough $T_1=T_1(M)\in(0,T)$ consider the unique solution $\lambda\in L^2([0,T_1])$ of the truncated fixed-point problem \eqref{Krylov_tr1}, \eqref{Krylov_tr2}. The corresponding solution $p$ of \eqref{Krylov_tr1} satisfies
\begin{equation}
\big\|\!-(\alpha+\lambda^{M,T_1}_t)\,\partial_y p\big\|_{L^2([0,T_1]\times[0,\infty))}\le(\alpha+M)\,C_2, 
\end{equation} 
where $C_2$ is as in \eqref{ux_bound}. Repeating the estimates from \eqref{short_time_heat_kernel} we get therefore
\begin{equation}
\begin{split}
\|(\partial_yp)(\cdot,0)\|_{L^2([0,T_1])}
\le (\alpha+M)\,C_2\,\bigg(\int_0^\infty \|(\partial_y\psi_\sigma)(\cdot,z,0)\|_{L^1([0,T_1])}^2\,\mathrm{d}z\bigg)^{1/2} \\
+\bigg\|\int_0^\infty f_\nu(z)\,(\partial_y\psi_\sigma)(t,z,0)\,\mathrm{d}z\bigg\|_{L^2([0,T_1])}.
\end{split}
\end{equation}
The fixed-point constraint \eqref{Krylov_tr2}, the lower bound in \eqref{eq:exp_est}, and the latter inequality give $\|\lambda\|_{L^2([0,T_1])}\le M$ upon decreasing the value of $T_1=T_1(M)\in(0,T)$ if necessary. Such a $T_1$ belongs to the set on the right-hand side of \eqref{max_exist_int}, since $\lambda^{M,T_1}=\lambda$ and consequently $\lambda$ is a solution of the fixed-point problem \eqref{FPprob1}, \eqref{FPprob2}.

\medskip

We show next that for every element $T_1$ of the set on the right-hand side of \eqref{max_exist_int} the corresponding solution of the fixed-point problem \eqref{FPprob1}, \eqref{FPprob2} is unique. To this end, 
for any two solutions $\lambda,\,\widetilde{\lambda}\in L^2([0,T_1])$ we let $M=1+\|\lambda\|_{L^2([0,T_1])}\vee\|\widetilde{\lambda}\|_{L^2([0,T_1])}$. Then, for any $\varepsilon\in(0,T_1]$ the restrictions of both $\lambda$ and $\widetilde{\lambda}$ to $[0,\varepsilon]$ solve the truncated fixed-point problem \eqref{Krylov_tr1}, \eqref{Krylov_tr2} on $[0,\varepsilon]$. Combining this observation with the contraction property established in Steps 1-3 we find an $\varepsilon\in(0,T_1]$ such that $\lambda_t=\widetilde{\lambda}_t$ for almost every $t\in[0,\varepsilon]$.  

\medskip

With this $\varepsilon$ and the solution $p\in W^{1,2}_2([0,\varepsilon]\times[0,\infty))$ of the Cauchy-Dirichlet problem in \eqref{Krylov_tr1} we consider the mapping which takes $L^2([\varepsilon,(2\varepsilon)\wedge T_1])$ functions $\rho$ to 
\begin{equation*}
-C\,\frac{\sigma^2}{2}\,\frac{(\partial_y u)(\cdot,0)}{\int_0^\infty u(\cdot,y)\,\mathrm{d}y}, 
\end{equation*}
where  $u$ is the unique solution of 
\begin{equation}
\begin{split}
\partial_t u=-(\alpha+\rho^{M,\varepsilon,(2\varepsilon)\wedge T_1}_t)\,\partial_y u+\frac{\sigma^2}{2}\,\partial^2_y u,\;\; u(\varepsilon,\cdot)=p(\varepsilon,\cdot),\;\; u(\cdot,0)=0,\\ 
u\in W^{1,2}_2([\varepsilon,(2\varepsilon)\wedge T_1]\times[0,\infty))
\end{split}
\end{equation}
and 
\begin{equation}
\begin{split}
\rho^{M,\varepsilon,(2\varepsilon)\wedge T_1}:=& \rho\,\mathbf{1}_{\{\|\rho\|_{L^2([\varepsilon,(2\varepsilon)\wedge T_1])}
\le(M^2-\|\lambda\|^2_{L^2([0,\varepsilon])})^{1/2}\}} \\
& +\rho\,\frac{(M^2-\|\lambda\|^2_{L^2([0,\varepsilon])})^{1/2}}{\|\rho\|_{L^2([\varepsilon,(2\varepsilon)\wedge T_1])}}\,\mathbf{1}_{\{\|\rho\|_{L^2([\varepsilon,(2\varepsilon)\wedge T_1])}
>(M^2-\|\lambda\|^2_{L^2([0,\varepsilon])})^{1/2}\}}.
\end{split}
\end{equation}
This mapping is well-defined with range contained in $L^2([\varepsilon,(2\varepsilon)\wedge T_1])$, since one can regard $u$ as the restriction of the unique solution of 
\begin{equation}\label{Krylov_PDE_piecewise}
\begin{split}
\partial_t u = -(\alpha+\xi^{M,(2\varepsilon)\wedge T_1}_t)\,\partial_y u +\frac{\sigma^2}{2}\,\partial_y^2 u,\;\;
u(0,\cdot)=f_\nu,\;\;u(\cdot,0)=0,\\
u\in W^{1,2}_2([0,(2\varepsilon)\wedge T_1]\times[0,\infty))
\end{split}
\end{equation}
to $[\varepsilon,(2\varepsilon)\wedge T_1]\times[0,\infty)$, where
\begin{equation}
\xi^{M,(2\varepsilon)\wedge T_1}_t:=\begin{cases}
\lambda_t & \mathrm{if}\;\;t\in[0,\varepsilon), \\
\rho^{M,\varepsilon,(2\varepsilon)\wedge T_1}_t & \mathrm{if}\;\;t\in[\varepsilon,(2\varepsilon)\wedge T_1],
\end{cases}
\end{equation}
and use the assumptions $f_\nu\in W^1_2([0,\infty))$, $f_\nu(0)=0$, the existence and uniqueness result of \cite[Chapter III, Remark 6.3]{LSU} for \eqref{Krylov_PDE_piecewise}, the well-definedness of the second evaluation map in \cite[Chapter II, Lemma 3.4]{LSU}, and the lower bound in \eqref{eq:exp_est}.

\medskip

Moreover, the described mapping is a contraction on $L^2([\varepsilon,(2\varepsilon)\wedge T_1])$. Indeed, repeating the analysis of Steps 1-3, replacing every occurence of the interval $[0,\varepsilon]$ by $[\varepsilon,(2\varepsilon)\wedge T_1]$ and estimating 
\begin{equation*}
\mathrm{ess\,sup}_{t\in [\varepsilon,(2\varepsilon)\wedge T_1]} \|(\partial_y u)(t,\cdot)\|_{L^2([0,\infty))},\;\;\|(\partial_y u)(\cdot,0)\|_{L^2([\varepsilon,(2\varepsilon)\wedge T_1])}
\end{equation*}
via the boundedness of the evaluation maps in \cite[Chapter II, Lemma 3.4]{LSU} and the results of \cite[Chapter III, Section 6]{LSU} for the problem \eqref{Krylov_PDE_piecewise} we conclude that the Lipschitz constant of the mapping does not exceed 
\begin{equation*}
C\,\frac{\sigma^2}{2}\,C_7\,\bigg(\bigg(\int_0^\infty \|(\partial_y\psi_\sigma)(\cdot,z,0)\|_{L^1([0,\varepsilon])}^2\,\mathrm{d}z\bigg)^{1/2}+\varepsilon^{1/2}\bigg),
\end{equation*}
where the constant $C_7$ is the same as in \eqref{eq:contraction}. It follows that $\lambda_t=\widetilde{\lambda}_t$ for almost every $t\in[\varepsilon,(2\varepsilon)\wedge T_1]$, as the restrictions of $\lambda$ and $\widetilde{\lambda}$ to $[\varepsilon,(2\varepsilon)\wedge T_1]$ are both fixed-points of the mapping in consideration. A sequential repetition of the same argument on the time intervals 
\begin{equation*}
[(2\varepsilon)\wedge T_1,(3\varepsilon)\wedge T_1],\;[(3\varepsilon)\wedge T_1,(4\varepsilon)\wedge T_1],\;\ldots
\end{equation*}
yields $\lambda_t=\widetilde{\lambda}_t$ for almost every $t\in[0,T_1]$.

\medskip

Part (b) of the proposition is an immediate consequence of the just established uniqueness assertion. In addition, the latter allows to combine the solutions of the fixed-point problem \eqref{FPprob1}, \eqref{FPprob2} for different elements $T_1$ of the set on the right-hand side of \eqref{max_exist_int} to a function $\lambda\in L^2_{loc}([0,t_{reg}))$, with $t_{reg}$ defined via \eqref{max_exist_int}. To obtain part (a) of the proposition it remains to check $\lim_{T_1\uparrow t_{reg}} \|\lambda\|_{L^2([0,T_1])}=\infty$ if $t_{reg}<T$.  If $t_{reg}<T$ and $\lim_{T_1\uparrow t_{reg}} \|\lambda\|_{L^2([0,T_1])}<\infty$ were to hold, then $\lambda\in L^2([0,t_{reg}])$ would be a solution of the fixed-point problem \eqref{FPprob1}, \eqref{FPprob2} on $[0,t_{reg}]$. In addition, with $p\in W^{1,2}_2([0,t_{reg}]\times[0,\infty))$ being the corresponding solution of the Cauchy-Dirichlet problem \eqref{FPprob1} the same arguments as in Steps 1-3 and the first paragraph of Step 4 would give the existence of a solution $\rho\in L^2([t_{reg},\widehat{T}])$ to the fixed-point problem
\begin{eqnarray}
&& \partial_t u = -(\alpha+\rho_t)\,\partial_y u + \frac{\sigma^2}{2}\,\partial_y^2 u,\;\; u(t_{reg},\cdot)=p(t_{reg},\cdot),\;\; u(\cdot,0)=0, \\
&& -C\,\frac{\sigma^2}{2}\,\frac{(\partial_y u)(t,0)}{\int_0^\infty u(t,y)\,\mathrm{d}y}=\rho_t\;\;\text{for\;almost\;every}\;\;t\in[t_{reg},\widehat{T}]
\end{eqnarray}
for $(\widehat{T}-t_{reg})\in(0,T-t_{reg})$ small enough (note that $p(t_{reg},\cdot)\in W^1_2([0,\infty))$ with $p(t_{reg},0)=0$ thanks to the well-definedness of the first evaluation map in \cite[Chapter II, Lemma 3.4]{LSU}). The concatenation of $\lambda$ and $\rho$ would then be a solution of the fixed-point problem \eqref{FPprob1}, \eqref{FPprob2} on $[0,\widehat{T}]$, a contradiction to the definition of $t_{reg}$ in \eqref{max_exist_int}. \ep

\medskip

Given an initial condition $\overline{Z}_0\stackrel{d}{=}\nu$ as in Theorem \ref{main1}, we define
\begin{equation}\label{eq:stoch_phys_sol}
\begin{split}
&\overline{Z}_t=\overline{Z}_0+\alpha\,t+\int_0^t \lambda_s\,\mathrm{d}s+\sigma\,\overline{B}_t,\;\; t\in[0,t_{reg}\wedge T\wedge\overline{\tau}),\\
&\overline{\tau}=\inf\,\{t\in[0,T]:\,\overline{Z}_t=0\},\quad\overline{Z}_t=0,\;\;t\in[\overline{\tau},t_{reg}\wedge T),
\end{split}
\end{equation}
with the pair $(t_{reg},\lambda)$ of Proposition \ref{prop:ex&uniq}(a).
The next proposition establishes that, until the explosion of the weak derivative of the cumulative loss process in the $L^2$ norm, any physical solution satisfying the conditions in Theorem \ref{main1}(a) and stopped upon hitting $0$ must be given by $\overline{Z}$.

\begin{proposition}\label{prop:stoch_ex_uniq}
Let $\nu$ be as in Theorem \ref{main1}. Then, 
\begin{enumerate}[(a)]
\item the process $\overline{Z}$ defined by \eqref{eq:stoch_phys_sol} satisfies the fixed-point constraint
\begin{equation}\label{eq:FP_constraint}
\lambda_t=C\,\partial_t\log \PP(\overline{\tau}>t)\;\;\text{for almost every}\;\;t\in[0,t_{reg}\wedge T); 
\end{equation}
\item for any physical solution $\overline{Y}$ satisfying the conditions in Theorem \ref{main1}(a), the corresponding time $t_{reg}>0$ and the stopped process $\overline{Z}_t:=\overline{Y}_{t\wedge\overline{\tau}}$, $t\in[0,t_{reg}\wedge T)$ are given by \eqref{eq:stoch_phys_sol}.
\end{enumerate}
\end{proposition}

\smallskip

\noindent\textbf{Proof.} For any $T_1\in(0,t_{reg}\wedge T)$ the argument employed in the proof of Lemma \ref{lemma:exp_bd} shows that the densities $p(t,\cdot)$, $t\in[0,T_1]$ of the restrictions of the laws of $\overline{Z}_t$, $t\in[0,T_1]$ to $(0,\infty)$, respectively, form a $W^{1,2}_2([0,T_1]\times[0,\infty))$-solution of \eqref{eq:cdp}. Consequently, the identity \eqref{eq:exp_test} and the lower bound in \eqref{eq:exp_est} reveal the function $t\mapsto\log\int_0^\infty p(t,y)\,\mathrm{d}y$ as absolutely continuous on $[0,T_1]$ with
\begin{equation}\label{key_rel}
\partial_t\log\int_0^\infty p(t,y)\,\mathrm{d}y=-\frac{\sigma^2}{2}\,\frac{(\partial_y p)(t,0)}{\int_0^\infty p(t,y)\,\mathrm{d}y}
\;\;\text{for almost every}\;\;t\in [0,T_1].
\end{equation}  
By combining \eqref{FPprob2} with \eqref{key_rel} we arrive at \eqref{eq:FP_constraint}, that is, part (a) of the proposition.  

\medskip

Next, we let $\lambda$ be the weak derivative of the loss function of a physical solution $\overline{Y}$ as in part (b) of the proposition and $t_{reg}>0$ be the explosion time of $\lambda$ in the $L^2$ norm. We also fix a $T_1\in(0,t_{reg}\wedge T)$ and denote by $p(t,\cdot)$, $t\in[0,T_1]$ the densities of the restrictions of the laws of $\overline{Y}_{t\wedge\overline{\tau}}$, $t\in[0,T_1]$ to $(0,\infty)$, respectively. Then, both \eqref{eq:FP_constraint} and \eqref{key_rel} hold. Moreover, substituting the right-hand side of \eqref{key_rel} for $\partial_t\log \PP(\overline{\tau}>t)$ in \eqref{eq:FP_constraint} we get \eqref{FPprob2}. Now, it follows from Proposition \ref{prop:ex&uniq}(b) that the pair $(t_{reg},\lambda)$ is the one of Proposition \ref{prop:ex&uniq}(a). Part (b) of the proposition at hand readily follows. \ep   

\medskip

We conclude the section with the proof of Theorem \ref{main1}.

\medskip

\noindent\textbf{Proof of Theorem \ref{main1}.} Parts (a) and (b) of the theorem follow directly from parts (a) and (b) of Proposition \ref{prop:stoch_ex_uniq}, respectively. Moreover, for any $T_1\in(0,t_{reg})$ the argument used in the proof of Lemma \ref{lemma:exp_bd} identifies the densities $p(t,\cdot)$, $t\in[0,T_1]$ of the restrictions of the laws of $\overline{Y}_{t\wedge\overline{\tau}}$, $t\in[0,T_1]$ to $(0,\infty)$ with the $W^{1,2}_2([0,T_1]\times[0,\infty))$-solution of \eqref{eq:cdp}. Thus, the restriction $(\partial_y p)(\cdot,0)\in L^2_{loc}([0,t_{reg}))$ of the weak derivative $\partial_y p$ is well-defined due to the well-definedness of the second evaluation map in \cite[Chapter II, Lemma 3.4]{LSU}, and the characterization \eqref{expl_crit} follows from \eqref{FPprob2}. \ep

\section{A priori regularity of physical solutions}
\label{se:regularity}

We begin this section by stating some elementary properties of physical solutions. 

\begin{lemma}\label{le:phys.sol.1}
Let $\overline{Y}$ be a c\`adl\`ag process satisfying (\ref{eq:phys_sol.1}), with the associated $\Lambda$, $\overline{\tau}^0$ and $\overline{\tau}$. Then, for $t\in(0,\overline{\tau}^0)$,
\begin{enumerate}[(a)]
\item the associated loss function $\Lambda$ is non-increasing; 
\item the laws of $\overline{Y}_t\,\mathbf{1}_{\{\overline{\tau}\ge t\}}$ and $\overline{Y}_{t-}\,\mathbf{1}_{\{\overline{\tau}\ge t\}}$, restricted to $(0,\infty)$, possess densities; the latter are bounded by a constant independent of $t$ if the law of $\overline{Y}_0$ possesses a bounded density; 
\item $\Lambda_{t-}= C\,\log \pp(\overline{\tau}\ge t)$; 
\item $\PP(\overline{Y}_{t-}\leq0)>0$.
\end{enumerate}
\end{lemma}

\smallskip

\noindent\textbf{Proof.} Property (a) is immediate from the definition of $\Lambda$. To deduce property (b) we notice that, for all $0<a<b<\infty$,
\begin{eqnarray*}
&& \PP\big(\overline{Y}_t\,\mathbf{1}_{\{\overline{\tau}\ge t\}}\in(a,b)\big) 
\leq \PP\big(\overline{Y}_0+\sigma\,\overline{B}_t\in(a-\alpha\,t-\Lambda_t,b-\alpha\,t-\Lambda_t)\big), \\
&& \PP\big(\overline{Y}_{t-}\,\mathbf{1}_{\{\overline{\tau}\ge t\}}\in(a,b)\big) 
\leq \PP\big(\overline{Y}_0+\sigma\,\overline{B}_t\in(a-\alpha\,t-\Lambda_{t-},b-\alpha\,t-\Lambda_{t-})\big),
\end{eqnarray*}
and the two right-hand sides are bounded above by a constant times $(b-a)$. This constant can be chosen to be the same for all values of $t$ in an interval bounded away from zero, and it is uniform for all $t\geq0$ if $\overline{Y}_0$ possesses a bounded density. To obtain property (c) we let $s\in[0,t)$, $\varepsilon\in(0,1)$ and employ the chain of estimates 
\begin{equation*}
\begin{split}
& \PP(\overline{\tau}>s) - \PP(\overline{\tau}\ge t)
\le \PP\Big(\overline{Y}_s>0,\,\inf_{r\in[s,t)} \overline{Y}_r \leq0\Big) \\
& \leq \PP\big(\overline{Y}_s\in(0,\varepsilon)\big)
+ \PP\Big(\inf_{r\in[s,t)} \big(\alpha\,(r-s)+\sigma\,(\overline{B}_r-\overline{B}_s) + \Lambda_r - \Lambda_s \big) \leq -\varepsilon\Big).
\end{split}
\end{equation*}
In view of the existence of $\Lambda_{t-}$ and property (b), the limit $\varepsilon\downarrow0$ of the limit superior $s\uparrow t$ of the latter upper bound is $0$, and property (c) readily follows. Finally, property (d) is a consequence of 
\begin{equation*}
\overline{Y}_{t-}\le\overline{Y}_0+\alpha\,t+\sigma\,\overline{B}_t,\;\;t\in(0,\overline{\tau}^0),
\end{equation*}
which is in turn due to property (a). \ep

\medskip

Let us fix an arbitrary c\`adl\`ag process $\overline{Y}$ satisfying (\ref{eq:phys_sol.1}), with the associated $\Lambda$ and $\overline{\tau}^0$. Recall that $p(t,\cdot)$ denotes the density of the distribution of $\overline{Y}_{t-}\,\mathbf{1}_{\{\overline{\tau}\ge t\}}$ restricted to $(0,\infty)$.
In the rest of the section, we establish certain regularity properties of $\overline{Y}$, which, ultimately, allow us to conclude that $\Lambda$ is H\"older continuous, with a H\"older exponent strictly greater than $1/2$, on any interval on which the normalized density at zero vanishes. 

\medskip

We assume that the law of $\overline{Y}_0$ admits a bounded density and begin with an auxiliary construction.
For fixed $t\in[0,\overline{\tau}^0)$ and $\varepsilon\in(0,\infty)$, we consider the sequence of processes $\overline{Y}^n$, $n\in\nn$ defined recursively as follows:
\begin{eqnarray}
&& \overline{Y}^1_s = \overline{Y}_{t-} + (\alpha\,s+\sigma\,\widetilde{B}_{s})\,\bone_{\{\overline{\tau}\geq t\}},\;\; s\in[0,\varepsilon], \\
&&\overline{Y}^n_s = \overline{Y}_{t-}+(\alpha\,s+\sigma\,\widetilde{B}_s+L^{n-1})\,\bone_{\{\overline{\tau}\geq t\}},\;\;s\in[0,\varepsilon],\;\;n\ge2, \\
&& L^n
= C\log\,\PP\Big(\overline{\tau}\geq t,\,\inf_{s\in[0,\varepsilon]} \overline{Y}^n_s > 0 \Big) - \Lambda_{t-},
\;\; n\geq 1, \label{eq.PhysSolReg.Ln.def}
\end{eqnarray}
where $\widetilde{B}_s:=\overline{B}_{t+s}-\overline{B}_t$, $s\in[0,\varepsilon]$ and $\Lambda_{0-}:=0$. The latter logarithm is well-defined, since $t<\overline{\tau}^0$ and $\widetilde{B}$ is independent of $\overline{Y}_s$, $s\in[0,t]$. By Lemma \ref{le:phys.sol.1}(c), $\overline{Y}^2_s\leq \overline{Y}^1_s$ for all $s\in[0,\varepsilon]$ with probability one. Then, by induction, we conclude that the sequences $\overline{Y}^n_s$, $n\in\nn$ are non-increasing for all $s\in[0,\varepsilon]$ with probability one. 

\begin{lemma}\label{le:PhysSolReg.4}
Suppose that the law of $\overline{Y}_0$ possesses a bounded density. Then, the following hold for any $t\in[0,\overline{\tau}^0)$.
\begin{enumerate}[(a)]
\item If $p(t,\cdot)$ satisfies
\begin{equation}\label{kappa_cond1}
\lim_{\eta\downarrow 0}\,\mathrm{ess\,sup}_{y\in(0,\eta)}\,p(t,y)=0,
\end{equation}
then there is a constant $C_L<\infty$ depending only on $C,\,\sigma,\,\Lambda_{t-},\,\|p(t,\cdot)\|_{L^\infty([0,\infty))}$ such that 
\begin{equation}
|L^n| \le C_L\,\varepsilon^{1/2}
\end{equation}
for all $n\in\nn$ sufficiently large and all $\varepsilon\in(0,\infty)$ sufficiently small, where $L^n$ is defined by \eqref{eq.PhysSolReg.Ln.def}.

\item If $p(t,\cdot)$ satisfies
\begin{equation}
p(t,y) \leq \widehat{C}\,y^\gamma,\;\;y\in(0,\eta)
\end{equation}
with some constants $\widehat{C}<\infty$, $\gamma\in(0,1]$, and $\eta>0$, then there is a constant $C_L<\infty$ depending only on
$C,\,\sigma,\,\widehat{C},\,\gamma,\,\eta,\,\Lambda_{t-},\,\|p(t,\cdot)\|_{L^\infty([0,\infty))}$ such that
\begin{equation}
|L^n| \leq C_L\,\varepsilon^{(1+\gamma)/2}
\end{equation}
for all $n\in\nn$ sufficiently large and all $\varepsilon\in(0,\infty)$ sufficiently small.
\end{enumerate}
\end{lemma}

\smallskip

\noindent\textbf{Proof.} Let $\kappa(y):=\mathrm{ess\,sup}_{z\in(0,y)} \,p(t,z)$, $y\in(0,\infty)$ and note that in the setting of part (b) it holds $\kappa(y)\leq\widehat{C}\,y^{\gamma}$, $y\in(0,\infty)$, where we have increased the value of $\widehat{C}$ if necessary (recall Lemma \ref{le:phys.sol.1}(b)). 
We have the estimates 
\begin{eqnarray*}
&& 0\geq e^{\Lambda_{t-}/C}\big(e^{L^1/C}-1\big)
= \pp\Big(\overline{\tau}\geq t,\,\inf_{s\in[0,\varepsilon]} \overline{Y}^1_s > 0 \Big)-\pp(\overline{\tau}\ge t) \\
&& = \int_0^\infty \Big(\pp\Big(\inf_{s\in[0,\varepsilon]} (\alpha\,s + \sigma\,\widetilde{B}_s)>-y\Big) - 1 \Big)\,p(t,y)\,\mathrm{d}y \\
&& \geq -2 \int_0^\infty \Phi\bigg(\frac{|\alpha|\,\varepsilon - y}{\sigma\sqrt{\varepsilon}}\bigg)\,p(t,y)\,\mathrm{d}y 
\geq -2\sqrt{\varepsilon} \int_0^{\infty} \Phi\left(1 - y/\sigma\right)\,p(t,y\sqrt{\varepsilon})\,\mathrm{d}y \\
&& = -2 \sqrt{\varepsilon} \int_0^{\iota/\sqrt{\varepsilon}} \Phi\left(1 - y/\sigma\right)\,p(t,y\sqrt{\varepsilon})\,\mathrm{d}y
-2\sqrt{\varepsilon} \int_{\iota/\sqrt{\varepsilon}}^\infty \Phi\left(1 - y/\sigma\right)\,p(t,y\sqrt{\varepsilon})\,\mathrm{d}y \\
&& \geq -2\sqrt{\varepsilon}\,\bigg(\int_0^{\infty} \Phi\left(1 - y/\sigma\right)\,\kappa(y\sqrt{\varepsilon})\,\mathrm{d}y+\|p(t,\cdot)\|_{L^\infty([0,\infty))}\,e^{-\varepsilon^{-1/4}}\bigg) \\
&& =:-\sqrt{\varepsilon}\,C_0(\varepsilon)
\end{eqnarray*}
for all $\iota\in(0,1)$ and sufficiently small $\varepsilon\in(0,\iota^3)$. Here, as before, $\Phi$ stands for the standard Gaussian cumulative distribution function. It is clear from \eqref{kappa_cond1} that $C_0(\varepsilon)\to 0$ as $\varepsilon\downarrow0$, and we conclude 
\begin{equation}
0\geq L^1 \geq C\log\big(1-e^{-\Lambda_{t-}/C}\sqrt{\varepsilon}\,C_0(\varepsilon)\big) 
\geq -2C e^{-\Lambda_{t-}/C}\sqrt{\varepsilon}\,C_0(\varepsilon)
\end{equation}
for all sufficiently small $\varepsilon>0$. In the setting of part (b) we have the additional upper bound
\begin{equation}
C_0(\varepsilon) \leq 2\varepsilon^{\gamma/2} \bigg(\widehat{C} \int_0^{\infty} \Phi\left(1 - y/\sigma\right)\,y^{\gamma}\,\mathrm{d}y + \|p(t,\cdot)\|_{L^\infty([0,\infty))}\bigg)=: C_1\,\varepsilon^{\gamma/2}.
\end{equation}

\smallskip

For $n\geq 2$, we find
\begin{equation}\label{eq.PhysSolReg.4.eq1}
\begin{split}
& e^{\Lambda_{t-}/C}\big(e^{L^n/C}-1\big)
= \pp\Big(\overline{\tau}\geq t,\,\inf_{s\in[0,\varepsilon]} \overline{Y}^n_s > 0 \Big)-\pp(\overline{\tau}\ge t) \\
& = \int_0^\infty \Big(\PP\Big(\inf_{s\in[0,\varepsilon]} (\alpha\,s + \sigma\,\widetilde{B}_s) + L^{n-1}> -y\Big) - 1 \Big)\,p(t,y)\,\mathrm{d}y \\
& \geq 
 -\int_0^{-L^{n-1}} p(t,y)\,\mathrm{d}y
-2 \int_{-L^{n-1}}^{\infty} \Phi\left(\frac{|\alpha|\varepsilon - y - L^{n-1}}{\sigma\sqrt{\varepsilon}}\right)\,p(t,y)\,\mathrm{d}y \\
& \geq 
 -\int_0^{-L^{n-1}} p(t,y)\,\mathrm{d}y
-2\sqrt{\varepsilon} \int_0^{\infty} \Phi\left(1 - y/\sigma\right)\,p(t,y\sqrt{\varepsilon} -L^{n-1})\,\mathrm{d}y.
\end{split}
\end{equation}

\smallskip

Next, we choose an $\iota\in(0,1)$ small enough, so that for all sufficiently small $\varepsilon\in(0,\iota^3)$ it holds 
\begin{equation}\label{iota_small}
2Ce^{-\Lambda_{t-}/C}\sqrt{\varepsilon}\,C_0(\varepsilon)
\leq\frac{2C_6}{1-C_5}\sqrt{\varepsilon}<\iota
\end{equation} 
where
\begin{equation}
\begin{split}
& C_6:=\frac{C C_3}{1-C_2\iota\kappa(\iota)-C_3\sqrt{\varepsilon}},\;\;
C_5:= C_4\,\kappa(\iota)<1,\;\;C_4:= \frac{C C_2}{1 - C_2\iota\kappa(\iota) - C_3 \sqrt{\varepsilon}}, \\
& C_3:=2e^{-\Lambda_{t-}/C}\,\|p(t,\cdot)\|_{L^\infty([0,\infty))} \int_0^\infty \Phi(1 - y/\sigma)\,\mathrm{d}y,\;\; C_2:= e^{-\Lambda_{t-}/C}.
\end{split}
\end{equation}
In particular, \eqref{iota_small} implies
\begin{equation}
L^1\geq -2Ce^{-\Lambda_{t-}/C}\sqrt{\varepsilon}\,C_0(\varepsilon)
\geq-\frac{2C_6}{1-C_5}\sqrt{\varepsilon}>-\iota.
\end{equation}
Assuming that 
\begin{equation}
L^{n-1}\geq -\frac{2C_6}{1-C_5}\sqrt{\varepsilon}>-\iota
\end{equation}
for some $n\ge2$, the overall estimate in \eqref{eq.PhysSolReg.4.eq1} yields
\begin{equation}
e^{\Lambda_{t-}/C}\big(e^{L^n/C}-1\big)
\geq -L^{n-1}\kappa(-L^{n-1})
-2\sqrt{\varepsilon}\,\|p(t,\cdot)\|_{L^\infty([0,\infty))} \int_0^{\infty} \Phi\left(1 - y/\sigma\right)\,\mathrm{d}y,
\end{equation}
so that
\begin{equation}
\begin{split}
& L^n \geq C\log\big( 1 - C_2 (-L^{n-1}) \kappa(-L^{n-1})
- C_3\sqrt{\varepsilon}\big)
\geq C_5 L^{n-1} - C_6 \sqrt{\varepsilon} \\
& \geq -\frac{2C_5C_6}{1-C_5}\sqrt{\varepsilon} - C_6 \sqrt{\varepsilon}
= -\frac{(1+C_5)C_6}{1-C_5}\sqrt{\varepsilon} \geq -\frac{2C_6}{1-C_5} \sqrt{\varepsilon}.
\end{split}
\end{equation}
Thus, by induction,
\begin{equation}
L^n \geq -\frac{2C_6}{1-C_5}\sqrt{\varepsilon}=:-C_7 \sqrt{\varepsilon} > -\iota,\;\;n\geq1.
\end{equation}

\smallskip

Finally, we apply (\ref{eq.PhysSolReg.4.eq1}) again to obtain for all sufficiently small $\varepsilon>0$:
\begin{equation}
\begin{split}
& e^{\Lambda_{t-}/C}\big(e^{L^n/C}-1\big) 
 \geq L^{n-1} \kappa(-L^{n-1})
-2\sqrt{\varepsilon}\int_0^{\iota/\sqrt{\varepsilon} - C_7} \Phi(1 - y/\sigma)\,\kappa((y+C_7)\sqrt{\varepsilon})\,\mathrm{d}y \\
&\qquad\qquad\qquad\qquad\quad -2\sqrt{\varepsilon}\,\|p(t,\cdot)\|_{L^\infty([0,\infty))} \int_{\iota/\sqrt{\varepsilon} - C_7}^\infty \Phi(1 - y/\sigma)\,\mathrm{d}y \\
&\qquad\qquad\qquad\quad\;\;\, \geq L^{n-1} \kappa(-L^{n-1}) - C_8\sqrt{\varepsilon}\,C_0(\varepsilon), 
\end{split}
\end{equation}
and, hence, $L^n \geq -C_5(-L^{n-1})-C_9\sqrt{\varepsilon}\,C_0(\varepsilon)$, $n\ge2$, for suitable constants $C_8,C_9<\infty$. Iterating the latter inequality we end up with 
\begin{equation}
0\geq L^n \geq - \frac{2C_9}{1-C_5}\sqrt{\varepsilon}\,C_0(\varepsilon)
\end{equation}
for all $n\in\nn$ sufficiently large. Both parts of the lemma readily follow. \ep

\medskip

Next, we use the sequence $L^n$, $n\in\nn$ to construct an auxiliary process $\widetilde{Y}$ admitting a comparison to the physical solution $\overline{Y}$.

\begin{lemma}\label{le:PhysSolReg.2}
Suppose that the law of $\overline{Y}_0$ possesses a bounded density and the assumptions of part (a) or part (b) of Lemma \ref{le:PhysSolReg.4} hold. Then, for all $\varepsilon\in(0,\infty)$ sufficiently small, there is a continuous process $\widetilde{Y}$ satisfying
\begin{equation}\label{eq.PhysSolReg.tildeY.def}
\widetilde{Y}_s=\overline{Y}_{t-} + (\alpha\,s +\sigma\,\widetilde{B}_{s} 
+\widetilde{L})\,\bone_{\{\overline{\tau}\geq t\}},\;\; s\in[0,\varepsilon],
\end{equation}
with
\begin{equation}\label{eq.PhysSolReg.tildeL.def}
\widetilde{L}
= C\log\,\PP\Big(\overline{\tau}\geq t,\,\inf_{s\in[0,\varepsilon]}\widetilde{Y}_u > 0\Big)- C\log\,\PP(\overline{\tau}\geq t).
\end{equation}
Moreover, 
\begin{equation}\label{eq.PhysSolReg.tildeL.Holder}
\widetilde{L} \geq - C_L\,\varepsilon^{(1+\gamma)/2},
\end{equation}
for all $\varepsilon\in(0,\infty)$ sufficiently small, where $C_L$ is as in the corresponding part of Lemma \ref{le:PhysSolReg.4} and $\gamma$ should be set to $0$ in the case of the setting of part (a) of Lemma \ref{le:PhysSolReg.4}.
\end{lemma}

\smallskip

\noindent\textbf{Proof.} By Lemma \ref{le:PhysSolReg.4}, for $\varepsilon\in(0,\infty)$ sufficiently small, the sequence $L^n$, $n\in\nn$ has a limit $\widetilde{L}$. 
Hence, the processes $\overline{Y}^n$, $n\in\nn$ converge uniformly on $[0,\varepsilon]$ to the process $\widetilde{Y}$ defined by \eqref{eq.PhysSolReg.tildeY.def} with probability one, so that $\inf_{s\in[0,\varepsilon]}\overline{Y}^n_s$, $n\in\nn$ tend almost surely to $\inf_{s\in[0,\varepsilon]} \widetilde{Y}_s$. Clearly, the conditional distribution of the latter random variable given $\{\overline{\tau}\ge t\}$ has no atoms and, hence, 
\begin{equation}
\lim_{n\rightarrow\infty}\,\PP\Big(\overline{\tau}\geq t,\,\inf_{s\in[0,\varepsilon]} Y^n_s > 0\Big)
= \PP\Big(\overline{\tau}\geq t,\,\inf_{s\in[0,\varepsilon]} \widetilde{Y}_s > 0\Big),
\end{equation}
which yields \eqref{eq.PhysSolReg.tildeL.def}. The estimate \eqref{eq.PhysSolReg.tildeL.Holder} follows directly from Lemma \ref{le:PhysSolReg.4}. \ep


\begin{lemma}\label{le:PhysSolReg.3}
Suppose that the law of $\overline{Y}_0$ possesses a bounded density and that $\Lambda$ is continuous at $t$ and on $[t,t+\varepsilon]$ for some $\varepsilon\in(0,\infty)$ as in Lemma \ref{le:PhysSolReg.2}. Then, with any solution $\widetilde{Y}$ of \eqref{eq.PhysSolReg.tildeY.def}, \eqref{eq.PhysSolReg.tildeL.def} for that $\varepsilon$, it holds 
\begin{equation}
\Lambda_{t+s} - \Lambda_{t} \geq \widetilde{L},\;\;
s\in[0,\varepsilon]\cap[0,\overline{\tau}^0-t).
\end{equation}
\end{lemma}

\smallskip

\noindent\textbf{Proof.} Suppose that there exists an $s\in[0,\varepsilon]\cap[0,\overline{\tau}^0-t)$ such that $\Lambda_{t+s}-\Lambda_t<\widetilde{L}$. Since $\widetilde{L}<0$, we must have $s>0$.
Due to the continuity of $\Lambda$ we can further find an $s'>0$ such that $\Lambda_{t+s'} - \Lambda_t = \widetilde{L}$ and $\Lambda_{t+s''}- \Lambda_t > \widetilde{L}$ for all $s''\in[0,s')$. Therefore, for any $s''\in[0,s']$, the definitions of $\overline{Y}$, $\widetilde{Y}$ and the properties of Brownian motion give
\begin{equation*}
\mathbf{1}_{\{\overline{\tau}>t+s''\}} - \mathbf{1}_{\{\overline{\tau}\geq t,\,\inf_{r\in[0,\varepsilon]}\widetilde{Y}_r > 0\}} \geq 0,
\quad \PP\Big(\mathbf{1}_{\{\overline{\tau}>t+s''\}} - \mathbf{1}_{\{\overline{\tau}\geq t,\,\inf_{r\in[0,\varepsilon]} \widetilde{Y}_r > 0\}} >0 \Big)>0.
\end{equation*}
Taking $s''=s'$ we end up with $\Lambda_{t+s'} - \Lambda_t > \widetilde{L}$, which is the desired contradiction. \ep

\medskip

The following proposition shows that the conditions of Lemma \ref{le:PhysSolReg.4} imply the H\"older continuity of the cumulative loss process. 

\begin{proposition}\label{prop:HolderL.givenHolderp}
Suppose that the law of $\overline{Y}_0$ possesses a bounded density and that for some $t_0\in(0,\overline{\tau}^0)$, $\Lambda$ is continuous on $[0,t_0)$ and the assumption of part (a) or part (b) of Lemma \ref{le:PhysSolReg.4} applies for all $t\in[0,t_0)$. 
Then, there exist $\widetilde{C}<\infty$ and $\gamma\in[0,1]$ such that
\begin{equation}\label{eq.PhysSolReg.L.Holder}
|\Lambda_t-\Lambda_s| \leq \widetilde{C}\,|t-s|^{(1+\gamma)/2},\;\;s,t\in[0,t_0),
\end{equation}
where $\gamma$ can be chosen strictly positive in the case of part (b) of Lemma \ref{le:PhysSolReg.4}.
\end{proposition}

\smallskip

\noindent\textbf{Proof.} Combining Lemmas \ref{le:PhysSolReg.4}, \ref{le:PhysSolReg.2}, \ref{le:PhysSolReg.3} we conclude that for any $t\in[0,t_0)$ there is a constant $C_L<\infty$ such that
\begin{equation}
0\geq \Lambda_{s}-\Lambda_{t} \geq -C_L(s-t)^{(1+\gamma)/2}
\end{equation}
holds for all $s$ in a right neighborhood of $t$. The proposition now follows by noting that the size of such neighborhoods can be chosen uniformly in $t$. \ep

\medskip

Finally, we recall the definition of $r^*_t$ from \eqref{r_star_t} and connect this quantity to the assumption in part (b) of Lemma \ref{le:PhysSolReg.4}. 

\begin{proposition}\label{prop:Holderp.givenZerop}
Suppose that $\overline{Y}_0$ has a bounded density vanishing in a neighborhood of $0$. Consider any $t_0\in(0,\overline{\tau}^0)$ for which $r^*_{t_0}=0$. Then, there exist $\widehat{C},\eta\in(0,\infty)$ and $\gamma\in(0,1]$ such that 
\begin{equation}
p(t,y)\leq\widehat{C}\,y^\gamma,\;\;y\in(0,\eta),\;\;t\in(0,t_0).
\end{equation}
\end{proposition}

\smallskip

\noindent\textbf{Proof.} The assumption $r^*_{t_0}=0$ and Definition \ref{def_phys_sol.new} of a physical solution imply that $\Lambda$ is continuous on $[0,t_0)$. Moreover, we note that the conditions in Lemma \ref{le:PhysSolReg.4}(a) are satisfied for all $t\in[0,t_0)$. Next, we fix $s\in(0,t_0)$ and $\chi,y_0\in(0,\infty)$, pick a function $\phi^\chi\in C^\infty([0,\infty))$ with values in $[0,1]$ and support contained in $(y_0,y_0+\chi)$, and define
\begin{equation}
g(t,y)=\E\big[\phi^\chi(\overline{Y}_{t\wedge\overline{\tau}} - y)\big],\;\;(t,y)\in[0,s]\times\rr.
\end{equation}
In addition, we let 
\begin{equation}
Z_t:=-\alpha\,t + \sigma\,(\overline{B}_{s-t} - \overline{B}_{s}) + \Lambda_{s-t} - \Lambda_{s},\;\; t\in[0,s],
\end{equation}
and write $\mathbb{F}^Z=(\mathcal{F}^Z_t)_{t\in[0,s]}$ for the filtration generated by $Z$.
We also consider the stopping time with respect to $\mathbb{F}^Z$:
\begin{equation}
\theta:=\big(\inf\{t\in[0,s]\,:\,Z_t\notin(-y_0,\chi_1\eta)\}\big)\wedge (\chi_2\eta^2)\wedge s,
\end{equation}
where $\chi_1,\chi_2,\eta\in(0,\infty)$ are constants to be specified below.
Note that, whenever $\overline{Y}_{s-t}\geq0$, we have
\begin{equation}
Z_t\geq \overline{Y}_{s-t} - \overline{Y}_{s},
\end{equation}
and whenever $\overline{Y}_{s-t},\,\overline{Y}_s\geq0$, we have
\begin{equation}
Z_t=\overline{Y}_{s-t} - \overline{Y}_{s}.
\end{equation}
The latter always holds on $\{t\leq\theta\}\cap\{\overline{Y}_s\geq y_0\}$.

\medskip

We claim that 
\begin{equation}
g(s-t\wedge\theta,Z_{t\wedge\theta}),\;\; t\in[0,s]
\end{equation}
is a martingale with respect to $\mathbb{F}^Z$. To this end, we use that, for any $t\in[0,s]$, $\overline{Y}_{s-t}$ is independent of $\mathcal{F}^Z_t$ and that
\begin{equation}
\bone_{\{\overline{Y}_s\geq y_0\}}\,\bone_{\{\inf_{r\in[0,s-t]}\overline{Y}_r>0\}}\,\bone_{\{t\leq\theta\}}
= \bone_{\{\overline{Y}_s\geq y_0\}}\,\bone_{\{\inf_{r\in[0,s]}\overline{Y}_r>0\}}\,\bone_{\{t\leq\theta\}},
\end{equation}
which yields
\begin{equation}\label{mart_expl}
\begin{split}
& \;\E\big[\phi^\chi(\overline{Y}_s)\,\bone_{\{\overline{Y}_s\geq y_0\}}\,\bone_{\{\inf_{r\in[0,s]}\overline{Y}_r>0\}}\,\big|\,\mathcal{F}^Z_t\big]\,\bone_{\{t\leq\theta\}} \\
& = \EE\big[\phi^\chi(\overline{Y}_{s-t} - Z_t)\,\bone_{\{\inf_{r\in[0,s-t]}\overline{Y}_r>0\}}\,\big|\,\mathcal{F}^Z_t\big]\,\bone_{\{t\leq\theta\}} \\
& = \EE\big[\phi^\chi(\overline{Y}_{s-t} - y)\,\bone_{\{\inf_{r\in[0,s-t]}\overline{Y}_r>0\}}\big]\big|_{y=Z_t}\,\bone_{\{t\leq\theta\}}
= g(s-t,Z_t)\,\bone_{\{t\leq\theta\}}, 
\end{split}
\end{equation}
where we have relied on $Z_t\geq -y_0$, $t\leq\theta$. The first expression in \eqref{mart_expl} is a martingale multiplied by $\bone_{\{t\leq\theta\}}$, so that the last expression in \eqref{mart_expl} stopped at $\theta$ is a martingale.


\medskip

Applying the optional sampling theorem we obtain
\begin{equation}
\begin{split}
& g(s,0)
= \EE\big[g(s-\theta, Z_{\theta})\big]  \leq \PP\big(Z_{\theta}\neq -y_0,\, \theta < s \big)\,\sup_{(t,z)\in [0,s]\times [-y_0,\chi_1\eta]}\, g(t,z) \\
&\qquad\qquad\qquad\qquad\qquad\quad\;\; + \EE\big[g(s-\theta,-y_0)\,\bone_{\{Z_{\theta}= -y_0\}}\big]
+ \EE\big[g(0,Z_s)\,\bone_{\{Z_{\theta} \neq -y_0,\,\theta = s\}}\big]\\
&\,\leq\, \PP\big(Z_\theta\neq -y_0,\, \theta<s \big) \,\sup_{(t,z)\in [0,s]\times [-y_0,\chi_1\eta]}\,g(t,z)
+ \sup_{t\in[0,s]}\,g(t,-y_0)
+ \sup_{z\in [-y_0,\chi_1\eta]}\,g(0,z),
\end{split}
\end{equation}
which implies further
\begin{equation}\label{chi1chi2_0}
\begin{split}
& \bigg(\int_{y_0}^{y_0+\chi} \phi^\chi(y)\,p(s,y)\,\mathrm{d}y\bigg)\\
& \leq\, \PP\big(Z_\theta\neq -y_0,\,\theta<s\big) \,\Big(\sup_{t\in [0,s]}\,\mathrm{ess\,sup}_{z\in[0,\chi_1\eta + y_0 + \chi]}\, p(t,z)\Big) \|\phi^\chi\|_{L^1([0,\infty))} \\
&\quad\, +\!\Big(\sup_{t\in [0,s]}\,\mathrm{ess\,sup}_{z\in[0,\chi]}\, p(t,z)\!\Big) \|\phi^\chi\|_{L^1([0,\infty))}
\!+\!\big(\mathrm{ess\,sup}_{z\in [0,\chi_1\eta+y_0+\chi]}\,p(0,z)\big)
\|\phi^\chi\|_{L^1([0,\infty))}.
\end{split}
\end{equation}

\smallskip

Letting $y_0\in(0,\eta)$ with $\eta\in(0,\infty)$ sufficiently small we get in the case $s\ge\chi_2\eta^2$:
\begin{equation}\label{chi1chi2_1}
\begin{split}
& \PP\big(Z_\theta \neq -y_0,\,\theta<s\big)
 \leq \PP\Big(\inf_{t\in[0,\chi_2\eta^2]} \big(\Lambda_{s-t} - \Lambda_s -\alpha\,t + \sigma\,(\overline{B}_{s-t}-\overline{B}_s)\big) > -\eta \Big) \\
&\qquad\qquad\qquad\qquad\quad\;\;
+\PP\Big(\sup_{t\in[0,\chi_2\eta^2]} \big(\Lambda_{s-t} - \Lambda_s - \alpha\,t + \sigma\,(\overline{B}_{s-t}-\overline{B}_s)\big) > \chi_1\eta \Big)  \leq \frac{1}{2},
\end{split}
\end{equation}
once we make $\chi_1,\chi_2\in(0,\infty)$ sufficiently large, uniformly in $\eta$ and $s\in(0,t_0)$. Hereby, the second inequality in \eqref{chi1chi2_1} relies on the $1/2$-H\"{o}lder continuity of $\Lambda$, which is in turn due to Proposition \ref{prop:HolderL.givenHolderp}. In the case $s<\chi_2\eta^2$, we obtain similarly
\begin{equation}\label{chi1chi2_2}
\begin{split}
\PP\big(Z_\theta\neq -y_0,\,\theta<s\big)
& \leq \PP\Big(\sup_{t\in[0,\chi_2\eta^2]} \big(\Lambda_{s-t} - \Lambda_s - \alpha\,t + \sigma\,(\overline{B}_{s-t}-\overline{B}_s)\big) > \chi_1\eta \Big) \leq \frac{1}{2}.
\end{split}
\end{equation}
Combining \eqref{chi1chi2_0}, \eqref{chi1chi2_1}, \eqref{chi1chi2_2}, and $y_0\in(0,\eta)$ we end up with 
\begin{equation}\label{chi1chi2_3}
\begin{split}
\|\phi^\chi\|_{L^1([0,\infty))}^{-1}
\int_{y_0}^{y_0+\chi} \phi^\chi(y)\,p(s,y)\,\mathrm{d}y
\leq\frac{1}{2}\,\Big(\sup_{t\in [0,s]}\,\mathrm{ess\,sup}_{z\in [0,(\chi_1+1)\eta+\chi]}\,p(t,z)\Big) \\
 + \Big(\sup_{t\in [0,s]}\,\mathrm{ess\,sup}_{z\in[0,\chi]}\,p(t,z)\Big)
+ \big(\mathrm{ess\,sup}_{z\in [0,(\chi_1+1)\eta+\chi]}\,p(0,z) \big).
\end{split}
\end{equation}

\smallskip

To finish the proof we consider $\chi\in(0,\eta)$ in \eqref{chi1chi2_3} and use $r^*_s=0$ to conclude 
\begin{equation}
\begin{split}
\mathrm{ess\,sup}_{z\in[0,\eta]}\,p(s,z)
\leq & \;\frac{1}{2} \,\Big(\sup_{t\in [0,s]}\,\mathrm{ess\,sup}_{z\in[0,(\chi_1+2)\eta]}\,p(t,z)\Big) \\
& \;+ \big(\mathrm{ess\,sup}_{z\in[0,(\chi_1+2)\eta]}\, p(0,z) \big).
\end{split}
\end{equation}
Replacing $s$ by $t\in[0,s]$ and taking the supremum over $t\in[0,s]$ we find therefore
\begin{equation}
\begin{split}
\sup_{t\in [0,s]}\,\mathrm{ess\,sup}_{z\in[0,\eta]}\,p(s,z)
\leq&\;\frac{1}{2} \,\Big(\sup_{t\in [0,s]}\,\mathrm{ess\,sup}_{z\in[0,(\chi_1+2)\eta]}\,p(t,z)\Big) \\
&\; + \big(\mathrm{ess\,sup}_{z\in[0,(\chi_1+2)\eta]}\, p(0,z) \big).
\end{split}
\end{equation}
An iteration of this inequality yields
\begin{equation}
\begin{split}
\sup_{t\in [0,s]}\,\mathrm{ess\,sup}_{z\in[0,\eta]}\,p(s,z)
\leq &\; \frac{1}{2^n}\,\Big(\sup_{t\in [0,s]}\,\mathrm{ess\,sup}_{z\in[0,(\chi_1+2)^n\eta]}\,p(t,z)\Big) \\
&\;+2 \big(\mathrm{ess\,sup}_{z\in[0,(\chi_1+2)^n\eta]}\, p(0,z) \big)
\end{split}
\end{equation}
for all $n\in\nn$. It remains to choose $\widetilde{\eta}\in(0,\infty)$ such that $p(0,\cdot)$ vanishes on $[0,\widetilde{\eta}]$, let $\eta\in(0,\widetilde{\eta}/(\chi_1+2))$, and select $n$ as the integer part of $\log(\widetilde{\eta}/\eta)/\log(\chi_1+2)$ to deduce 
\begin{equation}
\begin{split}
\sup_{t\in [0,s]}\,\mathrm{ess\,sup}_{z\in[0,\eta]}\,p(s,z) 
& \leq\, 2^{-\log (\widetilde{\eta}/\eta)/\log(\chi_1+2)+1} \Big(\sup_{t\in [0,s]}\,\mathrm{ess\,sup}_{z\in [0,\widetilde{\eta}]}\, p(t,z) \Big) \\
&=2^{-\log\widetilde{\eta}/\log(\chi_1+2)+1}\,\Big(\sup_{t\in [0,s]}\,\mathrm{ess\,sup}_{z\in [0,\widetilde{\eta}]}\, p(t,z)\Big)\, \eta^{\log 2/\log \lambda}.
\end{split}
\end{equation}
The proposition follows by noting that the factor in front of $\eta^{\log 2/\log \lambda}$ can be bounded by a constant $\widehat{C}\in(0,\infty)$ independent of $s\in(0,t_0)$. \ep

\medskip

Combining Proposition \ref{prop:Holderp.givenZerop} with Proposition \ref{prop:HolderL.givenHolderp} we get Theorem \ref{thm:Holder.loss}.

\section{Jumps of the cumulative loss process} \label{se:sufficientCond.Jump}

In this section, we prove Theorem \ref{main2.2}. Let us fix any c\`adl\`ag process $\overline{Y}$ satisfying (\ref{eq:phys_sol.1}), with the associated $\Lambda$, $\overline{\tau}^0$ and $r^{**}$, and with $\overline{Y}_0$ admitting a density.
For fixed $t\in[0,\overline{\tau}^0)$ and $\varepsilon\in(0,\infty)$, we define the Brownian motion $\widetilde{B}_s:=\overline{B}_{t+s}-\overline{B}_t$, $s\in[0,\varepsilon]$, and consider the processes $\widehat{Y}$ and $\widehat{Y}^n$, $n\in\nn$ defined as follows:
\begin{equation}\label{eq.suffJump.old.Ln.def}
\begin{split}
& \widehat{Y}_s = (\overline{Y}_{t-} + \alpha\,s+\sigma\,\widetilde{B}_{s}+ \widehat{L}_s)\,\bone_{\{\overline{\tau}\geq t\}},\;\; s\in[0,\varepsilon], \\
& \widehat{L}_s = \Lambda_{t+s} - \Lambda_{t-} = C\log\,\PP\Big(\inf_{r\in[0,s]}\widehat{Y}_r > 0 \Big) - C\log\,\PP(\overline{\tau}\geq t), \;\; s\in[0,\varepsilon]\\
& \widehat{Y}^1_s = (\overline{Y}_{t-} + \alpha\,s+\sigma\,\widetilde{B}_{s})\,\bone_{\{\overline{\tau}\geq t\}},\;\; s\in[0,\varepsilon], \\
& \widehat{Y}^n_s = (\overline{Y}_{t-}+\alpha\,s+\sigma\,\widetilde{B}_s + \widehat{L}^{n-1})\,\bone_{\{\overline{\tau}\geq t\}},\;\;s\in[0,\varepsilon],\;\;n\ge2, \\
& \widehat{L}^n = C\log\,\PP(\widehat{Y}^n_{\varepsilon} > 0) - C\log\,\PP(\overline{\tau}\geq t),\;\; n\geq 1, 
\end{split}
\end{equation}
where $\Lambda_{0-}:=0$. It is clear that $\widehat{Y}^1_\varepsilon \geq \widehat{Y}_\varepsilon$. In addition, if $\widehat{Y}^n_\varepsilon\geq \widehat{Y}_\varepsilon$, then 
\begin{equation}
\mathbf{1}_{\{\widehat{Y}^{n}_\varepsilon>0\}} \geq \mathbf{1}_{\{\inf_{r\in[0,\varepsilon]} \widehat{Y}_r > 0\}},
\end{equation}
so that $\widehat{Y}^{n+1}_\varepsilon\geq \widehat{Y}_\varepsilon$. Thus, by induction $\widehat{Y}^n_\varepsilon\geq\widehat{Y}_\varepsilon$, $n\in\nn$, and
\begin{equation}\label{eq.SuffCondJump.hatL.est}
\widehat{L}_\varepsilon \leq \widehat{L}^n
= C\big(\log \PP(\widehat{Y}^n_{\varepsilon} > 0) - \log\PP(\overline{\tau}\geq t)\big),\;\;n\in\nn.
\end{equation}

\smallskip

Suppose that $r^{**}_t>c^*/C$ for some $c^*<\infty$ to be specified below, but $\Lambda$ is continuous at $t$.
Then, there exists an $\eta\in(0,1)$ such that
\begin{equation}
p(t,y) \geq \frac{c^*}{C}\,\PP(\overline{\tau}\geq t)=:\iota,\;\; y\in(0,\eta).
\end{equation}
It is easy to see that the distribution of $\widehat{Y}^1_\varepsilon$ restricted to $\RR\setminus\{0\}$ has a density, which we denote by $h$. Moreover, for $M\in(0,\infty)$ we pick an $\varepsilon_0=\varepsilon_0(M,\alpha,\eta)>0$ such that 
\begin{equation}
\big((M+1)\sqrt{\varepsilon},\, 2\eta/3 \big) \subset \big(M\sqrt{\varepsilon}+\alpha\,\varepsilon,\,\eta - M\sqrt{\varepsilon}+\alpha\,\varepsilon\big),\;\;\varepsilon\in(0,\varepsilon_0).
\end{equation}
Then, for any $\varepsilon\in(0,\varepsilon_0)$ and $y\in\big((M+1)\sqrt{\varepsilon},\,2\eta/3\big)$ we have
\begin{equation}\label{eq.SuffCondJump.g.est}
\begin{split}
h(y) & =  \frac{1}{\sqrt{2\pi}\sigma} \int_{(- y+\alpha\varepsilon)/\sqrt{\varepsilon}}^\infty e^{-z^2/(2\sigma^2)}\,p(t,z\sqrt{\varepsilon}+y-\alpha\varepsilon)\,\mathrm{d}z \\
& \geq \frac{\iota}{\sqrt{2\pi}\sigma} \int_{(- y + \alpha\varepsilon)/\sqrt{\varepsilon}}^{(\eta - y + \alpha\varepsilon)/\sqrt{\varepsilon}} e^{-z^2/(2\sigma^2)}\,\mathrm{d}z \\
& \geq  \frac{\iota}{\sqrt{2\pi}\sigma} \int_{-M}^M e^{-z^2/(2\sigma^2)}\,\mathrm{d}z
= \frac{c^*}{C}\,\PP(\overline{\tau}\geq t)\,\big(2\Phi(M/\sigma)-1\big),
\end{split}
\end{equation}
where $\Phi$ is the standard Gaussian cumulative distribution function. In addition,
\begin{equation}\label{Phi32}
\begin{split}
& \;\PP(\overline{\tau}\geq t) - \PP(\widehat{Y}^1_\varepsilon>0)
= \PP(\overline{\tau}\geq t,\,\widehat{Y}^1_\varepsilon\leq 0) \\
& \geq \frac{1}{\sqrt{2\pi}\sigma} \int_{-\infty}^0 \int_{(- y + \alpha\varepsilon)/\sqrt{\varepsilon}}^{(\eta - y + \alpha\varepsilon)/\sqrt{\varepsilon}} e^{-z^2/(2\sigma^2)}\,p(t,z\sqrt{\varepsilon} + y - \alpha\varepsilon)\,\mathrm{d}z\,\mathrm{d}y \\ 
& \geq \frac{\iota}{\sqrt{2\pi}\sigma} \int_{-\infty}^0 \int_{(- y + \alpha\varepsilon)/\sqrt{\varepsilon}}^{(\eta - y + \alpha\varepsilon)/\sqrt{\varepsilon}} e^{-z^2/(2\sigma^2)}\,\mathrm{d}z\,\mathrm{d}y \\
& \geq \frac{\iota\sqrt{\varepsilon}}{\sqrt{2\pi}\sigma} \int_{-1}^0 \int_{(\alpha\varepsilon)/\sqrt{\varepsilon} - y}^{(\eta + \alpha\varepsilon)/\sqrt{\varepsilon} - y} e^{-z^2/(2\sigma^2)}\,\mathrm{d}z\,\mathrm{d}y \\
& \geq \frac{\iota\sqrt{\varepsilon}}{\sqrt{2\pi}\sigma} \int_{\alpha\sqrt{\varepsilon} + 1}^{\eta/\sqrt{\varepsilon}+\alpha\sqrt{\varepsilon}} e^{-z^2/(2\sigma^2)}\,\mathrm{d}z
\geq \iota\sqrt{\varepsilon}\,\big(\Phi(3/\sigma) - \Phi(2/\sigma)\big),
\end{split}
\end{equation}
upon decreasing the value of $\varepsilon_0>0$ if necessary. 

\medskip

The estimate \eqref{Phi32} implies
\begin{equation}
\begin{split}
& \;q^1:= \log\PP(\overline{\tau}\geq t) 
- \log\PP(\widehat{Y}^1_\varepsilon>0) \\
& \geq  \frac{1}{\PP(\overline{\tau}\geq t)}\,\iota\sqrt{\varepsilon} \,\big(\Phi(3/\sigma) - \Phi(2/\sigma)\big) 
=\frac{c^*}{C}\sqrt{\varepsilon}\,\big(\Phi(3/\sigma) - \Phi(2/\sigma)\big).
\end{split}
\end{equation}
We now choose $M\in(0,\infty)$ and $c^*=c^*(\sigma)<\infty$ satisfying
\begin{equation}
c^*\big(\Phi(3/\sigma)-\Phi(2/\sigma)\big) > M+1,
\quad c^*\big(2\Phi(M/\sigma)-1\big) > 1.
\end{equation}
With such $M$, $c^*$ we see from \eqref{eq.SuffCondJump.g.est}:
\begin{equation}
h(y) \geq \frac{c^*}{C}\,\PP(\overline{\tau}\geq t)\,\big(2\Phi(M/\sigma)-1\big)=:\frac{C_1}{C}\,\PP(\overline{\tau}\geq t)
\end{equation}
for all $y\in(C q^1,2\eta/3)$, where $C_1>1$. 

\medskip

We let
\begin{equation}
q^n := \log \PP(\overline{\tau}\geq t) - \log \PP(\widehat{Y}^n_\varepsilon>0)
\end{equation}
and note that $\widehat{Y}^n_{\varepsilon} =\widehat{Y}^{1}_{\varepsilon} - C q^{n-1}$ on $\{\overline{\tau}\geq t\}$. Hence, the density $h^n$ of the distribution of $\widehat{Y}^n_{\varepsilon}\,\bone_{\{\overline{\tau}\geq t\}}$ restricted to $\RR\setminus\{0\}$ is given by 
\begin{equation}
h^n(y) = h(y + Cq^{n-1}),\;\;y\neq 0.
\end{equation}
Assuming further that $(q^k)_{k=1}^{n-1}$ is a non-decreasing sequence we obtain $\widehat{Y}^n_{\varepsilon} \leq \widehat{Y}^{n-1}_{\varepsilon}$ and
\begin{equation}
q^{n} - q^{n-1}
= \log\PP(\widehat{Y}^{n-1}_{\varepsilon}>0)
- \log \PP(\widehat{Y}^n_{\varepsilon}>0)\geq 0.
\end{equation}
By induction we conclude that $(q^n)_{n\ge 1}$ is non-decreasing and $(\widehat{Y}^n_\varepsilon)_{n\ge 1}$ is non-increasing. 

\medskip

Suppose that $q^n\leq 2\eta/3$ for all $n\geq1$. Then, for $n\geq 3$:
\begin{equation}
\begin{split}
& \;\PP(\widehat{Y}^{n-1}_\varepsilon>0)-\PP(\widehat{Y}^{n}_\varepsilon>0)
=\PP(\widehat{Y}^{n-1}_\varepsilon>0,\,\widehat{Y}^n_\varepsilon\le 0)
= \PP\big(0<\widehat{Y}^{n-1}_\varepsilon\leq C (q^{n-1}-q^{n-2})\big) \\
& = \int_{0}^{C (q^{n-1}-q^{n-2})} h^{n-1}(y)\,\mathrm{d}y
= \int_{C q^{n-2}}^{C q^{n-1}} h(y)\,\mathrm{d}y
\geq (q^{n-1} - q^{n-2})\,C_1\,\PP(\overline{\tau}\geq t),
\end{split}
\end{equation}
and, thus, $q^n - q^{n-1} \geq C_1\,(q^{n-1} - q^{n-2})$, $n\ge3$. Moreover, it is easy to see that $q^2-q^1>0$, which in view of $C_1>1$ yields a contradiction to $q^n\leq 2\eta/3$, $n\ge1$. Consequently, for all $\varepsilon_0=\varepsilon_0(M,\alpha,\eta)>0$ sufficiently small, one can find for each $\varepsilon\in(0,\varepsilon_0)$ an $n\in\nn$ such that $q^n > 2\eta/3$. At this point, \eqref{eq.SuffCondJump.hatL.est} gives 
\begin{equation}
-\widehat{L}_\varepsilon
\geq C \sup_{n\in\nn} \big(\log\PP(\overline{\tau}\geq t)
- \log \PP(\widehat{Y}^n_\varepsilon>0) \big)
= C \sup_{n\in\nn} q^n > 2C\eta/3.
\end{equation}
This is the desired contradiction to the continuity of $\Lambda$ at $t$.

\section{Convergence of the finite-particle systems}
\label{se:convergence}

This last section is devoted to the proof of Theorem \ref{main2}, and we work under the assumptions of that theorem throughout the rest of the paper. Recall the equations \eqref{dynamics_w_defaults.2.1}, \eqref{dynamics_w_defaults.2.2}, \eqref{dynamics.YN1}, satisfied by the finite-particle system $(\widetilde{Y}^1,\widetilde{Y}^2,\ldots,\widetilde{Y}^N)$ and the definition \eqref{eq.V.asCascadeSize.def} of the cascade sizes, needed to uniquely determine the finite-particle system. We note that the equations \eqref{dynamics_w_defaults.2.1}, \eqref{dynamics_w_defaults.2.2} guarantee that the processes $\widetilde{Y}^1,\,\widetilde{Y}^2,\,\ldots,\,\widetilde{Y}^N$ never jump across $-1$. Our proof follows the line of reasoning in \cite{Delarue2}, and we often refer to the results established therein.

\medskip

Our first aim is to establish the tightness of the sequence of empirical measures $\widetilde{\mu}^N=\frac{1}{N}\,\sum_{i=1}^N \delta_{\widetilde{Y}^i}$, $N\in\nn$. To this end, we start with the following lemma, which is the analogue of \cite[Lemma 5.2]{Delarue2}.

\begin{lemma}\label{le:E.Holder}
For any $\chi>0$, there exists some $\upsilon=\upsilon(\chi)\in(0,1)$ (independent of $N$) such that
\begin{equation}
\PP\big(\exists\,t\in[0,\upsilon]:\;1-S_t/N \geq \upsilon^{-1} t^{1/4}\big) \leq \chi.
\end{equation}
\end{lemma}

\smallskip

\noindent\textbf{Proof.} We introduce the auxiliary particle system $(\widehat{Y}^1,\widehat{Y}^2,\ldots,\widehat{Y}^N)$ defined analogously to $(\widetilde{Y}^1,\widetilde{Y}^2,\ldots,\widetilde{Y}^N)$, but with the equation \eqref{dynamics.YN1} replaced by 
\begin{equation}
\widehat{Y}_t^i=\widetilde{Y}_0^i+\alpha\,t+\sigma\,B^i_t 
-\frac{1+C}{2}\,(1-\widehat{S}_t/N),\;\; t\in[0,T].
\end{equation}
More specifically, we substitute \eqref{eq.V.asCascadeSize.def} by
\begin{equation}
\widehat{D}_t=\bigg(\inf\bigg\{k=1,\,2,\,\ldots,\,N-\widehat{S}_{t-}:\;
\widehat{Y}^{(k)}_{t-}-\frac{1+C}{2}\,\frac{k-1}{N} > 0\bigg\}-1\bigg)\wedge \widehat{S}_{t-},
\end{equation}
and rewrite \eqref{dynamics_w_defaults.2.1} for $(\widehat{Y}^1,\widehat{Y}^2,\ldots,\widehat{Y}^N)$ accordingly.
Fix an arbitrary $\chi_1\in(0,1)$. By repeating the proof of \cite[Lemma 5.2]{Delarue2} we conclude that for any $\chi>0$ there exists some $\upsilon=\upsilon(\chi)\in(0,1)$ (independent of $N$) such that
\begin{equation}\label{chi_est1}
\pp\big(\exists\,t\in[0,\chi_1^4\upsilon^4]:\;1-\widehat{S}_t/N\ge \upsilon^{-1}t^{1/4}\big)\le\chi. 
\end{equation}
In fact, in \cite{Delarue2} each particle shifts the locations of the other particles every time it hits a new integer, which makes $\widehat{S}$ even smaller. This observation allows to simplify some parts of the proof of \cite[Lemma 5.2]{Delarue2} when deriving the estimate \eqref{chi_est1}.

\medskip

If $\chi_1\in(0,1)$ is chosen for the following to hold:
\begin{equation}
-\frac{1+C}{2}\,y_1\leq C\log(1-y_1), \;\; 
-\frac{1+C}{2}\,y_1\leq C\log\bigg(1-\frac{y_1}{1-y_2}\bigg),\;\; y_1,y_2\in[0,\chi_1),
\end{equation}
then on the complement of the event in \eqref{chi_est1} we have for all $t\in[0,\chi_1^4\upsilon^4]$ and all $k=1,\,2,\,\ldots,\,\widehat{D}_t+1$:
\begin{equation}
-\frac{1+C}{2}\,(1-\widehat{S}_t/N) \leq C\log(\widehat{S}_t/N),\;\;
-\frac{1+C}{2}\,\frac{k-1}{N} \leq C\log\bigg(1-\frac{k-1}{\widehat{S}_{t-}}\bigg).
\end{equation}
These inequalities and induction along the hitting times of zero for the auxiliary particles yield on the complement of the event in \eqref{chi_est1} for all $t\in[0,\chi_1^4 \upsilon^4]$ and all $i=1,\,2,\,\ldots,\,N$:
\begin{equation}
\widetilde{Y}^i_{t\wedge\tau^i} \geq \widehat{Y}^i_{t\wedge\widehat{\tau}^i},\;\; D_t \leq \widehat{D}_t,\;\; S_t \geq \widehat{S}_t.
\end{equation}
The lemma follows upon decreasing the value of $\upsilon\in(0,1)$ if necessary. \ep

\medskip

The next lemma is needed to prove the upper bound on jump sizes in the definition of a physical solution. It is the analogue of \cite[Lemma 5.3]{Delarue2}, and its proof is a simplified version of the proof of the latter. The present setting allows for a simplification, because each particle can only contribute to the cumulative loss process once, that is, it can only ``spike" once, in the terminology of \cite{Delarue2}.

\begin{lemma}\label{le:jump.necessary}
There exist some $C_0<\infty$, $\varepsilon>0$ such that for all $r\in(0,1)$, $t\in[0,T)$ and $s\in(0,(T-t)\wedge\varepsilon)$ one can find an $N_0=N_0(r,s)\in\nn$ with 
\begin{equation}\label{eq:jump.necessary}
\begin{split}
& \PP\bigg(\!\frac{S_{t-}}{N}\!\ge\!r,\,\forall \iota\!\leq\!\bigg(\!\frac{S_{t-}\!-\!S_{t+s}}{S_{t-}}\!-\!2s^{1/4}\!\bigg)^{\!\!+}\!\!:
\frac{1}{S_{t-}}\big|\{\tau^i\!\ge\! t,\widetilde{Y}^i_{t-}\!+\! C\log(1\!-\!\iota)\!\le\! 2s^{1/4}\}\big|\!\ge\!\frac{\iota}{1\!+\!s^{1/4}}\!\bigg)  \\
& \geq \PP\bigg(\frac{S_{t-}}{N}\!\ge\!r\bigg) - C_0 s
\end{split}
\end{equation}
for all $N\ge N_0$.
\end{lemma}

\smallskip

\noindent\textbf{Proof.} For an $\varepsilon>0$ we consider arbitrary $r\in(0,1)$, $t\in[0,T)$, $s\in(0,(T-t)\wedge\varepsilon)$ and work throughout on the event $\{S_{t-}/N\ge r\}$ (in particular, all events are intersected with $\{S_{t-}/N\ge r\}$, and all complements are taken with respect to $\{S_{t-}/N\ge r\}$). 
Then, for any $k\in\{0,1,\ldots,S_{t-}-S_{t+s}\}$ we have $\sum_{i=1}^N \bone_{A^{i,1}(k)} \geq k$, where
\begin{equation}
A^{i,1}(k):=\bigg\{\tau^i\geq t,\;\widetilde{Y}^i_{t-} + C\log\bigg(1-\frac{k}{S_{t-}}\bigg) -\alpha^-s+\sigma\inf_{s'\in[0,s]} (B^i_{t+s'}-B^i_t) \leq 0 \bigg\}.
\end{equation}
In addition, we define the events
\begin{equation}
A=\bigg\{\frac{1}{S_{t-}}\,\sum_{i=1}^N \bone_{A^{i,2}} \leq s\bigg\},
\;\; A^{i,2} = \bigg\{\tau^i\geq t,\;-\alpha^-s+\sigma\inf_{s'\in[0,s]} (B^i_{t+s'}-B^i_t) < -s^{1/4} \bigg\},
\end{equation}
let $\iota\in\big[0,\big(\frac{S_{t-}-S_{t+s}}{S_{t-}}\,(1+2s^{1/4}) - 2s^{1/4}\big)^+\big]$, and choose $k$ as the integer part of $\frac{\iota+s^{1/4}}{1+s^{1/4}}\,S_{t-}$, so that $0\le k\le S_{t-}-S_{t+s}$. Moreover, on $A^{i,1}(k)\cap (A^{i,2})^c$:
\begin{equation}
\begin{split}
\widetilde{Y}^i_{t-} + C \log (1-\iota) & \leq \widetilde{Y}^i_{t-} 
+ C \log \bigg(1-\frac{k}{S_{t-}}\,(1+s^{1/4}) + s^{1/4}\bigg) \\ 
& \leq \widetilde{Y}^i_{t-} 
+ C \log \bigg(1-\frac{k}{S_{t-}}\bigg) +  s^{1/4} \leq 2s^{1/4}.
\end{split}
\end{equation}
Consequently, on $A$:
\begin{equation}
\begin{split}
\frac{1}{S_{t-}}\,\sum_{i=1}^{N} \bone_{\{\tau^i\geq t,\,\widetilde{Y}^i_{t-} + C\log (1-\iota)\leq 2s^{1/4}\}} 
& \geq \frac{1}{S_{t-}}\,\sum_{i=1}^{N} \bone_{A^{i,1}(k)} 
- \frac{1}{S_{t-}}\,\sum_{i=1}^{N} \bone_{A^{i,2}} \\
& \geq \frac{k}{S_{t-}} - s 
\geq \frac{\iota + s^{1/4}}{1+s^{1/4}} - \frac{1}{Nr} - s
\geq \frac{\iota}{1+s^{1/4}}
\end{split}
\end{equation}
for $N \geq r^{-1}s^{-1}$, provided $\varepsilon$ is smaller than an appropriate uniform constant. Finally, we obtain $\pp(A^c)\le C_0s$ for some $C_0<\infty$ by conditioning on the information up to time $t$ and using Markov's inequality in conjunction with a standard estimate for Brownian motion. \ep

\medskip

Next, we recall the space $\mathcal{D}([0,T+1])$ of c\`{a}dl\`{a}g functions on $[0,T+1]$ that are continuous at $T+1$, endowed with the Skorokhod M1 topology (see, e.g., \cite{Delarue2}, \cite{Skorohod}, \cite{JacodShiryaev}, \cite{Whitt} for a detailed discussion of the M1 topology). We write $\mathcal{P}(\mathcal{D}([0,T+1]))$ for the space of probability measures on $\mathcal{D}([0,T+1])$, endowed with the topology of weak convergence. 

\begin{proposition}\label{prop:tightness}
The sequence $\widetilde{\mu}^N$, $N\in\nn$ is tight on $\mathcal{P}(\mathcal{D}([0,T+1]))$.
\end{proposition}

\smallskip

\noindent\textbf{Proof.} We first claim that the sequence of $\widetilde{Y}^1$, indexed by $N\in\nn$, is tight on $\mathcal{D}([0,T+1])$. To this end, we decompose $\widetilde{Y}^1$ into the sum of its continuous and jump parts. Notice that, between the jump times, $\widetilde{Y}^1$ is given by
\begin{equation}
\widetilde{Y}^1_t = \frac{\widetilde{Y}^1_0 + \alpha\,t + \sigma\,B^1_t + \bone_{\{\widetilde{Y}^1_t\geq -1\}}\,C \log( 1 - S_t/N)}{ 1 - \bone_{\{-1\leq \widetilde{Y}^1_t\leq 0\}}\,C \log(1 - S_t/N)}.
\end{equation}
Hence, the modulus of continuity of the continuous part of $\widetilde{Y}^1$ is bounded above by the modulus of continuity of a Brownian motion with drift, started from $\widetilde{Y}^1_0$. The same is true for the supremum of the absolute value of the continuous part of $\widetilde{Y}^1$. Moreover, the jump part of $\widetilde{Y}^1$ is non-increasing, and the supremum of its absolute value is bounded above by the supremum of the continuous part plus $1$, since $\widetilde{Y}^1$ does not jump across $-1$. The tightness of the sequence of $\widetilde{Y}^1$, indexed by $N\in\nn$, can be now deduced as in the proof of \cite[Lemma 5.4]{Delarue2} by invoking our Lemma \ref{le:E.Holder}. To obtain the proposition from this, it remains to use a standard argument from the theory of propagation of chaos, cf. \cite[Proposition 2.2]{Sznitman}. \ep

\medskip

We proceed with the proof of Theorem \ref{main2}. Let us denote by $\omega$ the canonical process in $\mathcal{D}([0,T+1])$ and introduce
\begin{equation}
m_t := \bone_{\{\inf_{s\in[0,t]} \omega_s \leq 0\}},\;\; t\in[0,T+1].
\end{equation}
By Proposition \ref{prop:tightness} the sequence $\widetilde{\mu}^N$, $N\in\nn$ is tight, and we write $\Pi_\infty$ for the law of an arbitrary limit point. In the remainder of the proof we assume that all limits are taken along the corresponding convergent subsequence of $\widetilde{\mu}^N$, $N\in\nn$. By repeating the first part of the proof of \cite[Theorem 4.4]{Delarue2} we find a countable set $J\subset [0,T+1]$ such that for any $t\in J^c$ it holds $\langle \mu,m_{t-} \rangle=\langle \mu,m_{t} \rangle$ and $\mu(\omega_{t-}=\omega_t)=1$ for $\Pi_{\infty}$-almost every $\mu\in\mathcal{P}(\mathcal{D}([0,T+1]))$. Hereby, $\langle \cdot,\cdot\rangle$ stands for the integral of a function of the canonical process with respect to a given measure. The following lemma is the analogue of \cite[Lemma 5.9]{Delarue2}, and its proof is postponed to Subsection \ref{subse:cont.mum}.

\begin{lemma}\label{le:cont.mum}
For $\Pi_{\infty}$-almost every $\mu\in\mathcal{P}(\mathcal{D}([0,T+1]))$ and any $\mu_n$, $n\in\nn$ converging weakly to $\mu$ we have
\begin{equation}
\lim_{n\rightarrow\infty}\,\langle \mu_n,m_t \rangle = \langle \mu,m_t \rangle,\;\;t\in J^c.
\end{equation}
\end{lemma}

Next, we fix a rational $T'\in J^c\cap[0,T]$, an integer $\ell\ge1$, elements $0=s_0<s_1<\cdots<s_\ell<T'$ of $J^c$, and uniformly continuous bounded functions $g_1,g_2,\ldots,g_\ell:\,\rr\to\rr$. In addition, for any uniformly continuous bounded function $G:\,\RR\rightarrow\RR$ we let
\begin{equation}
Q^N:= \EE\bigg[\widetilde{G}\bigg(\!\bigg\langle \widetilde{\mu}^N,\,
\prod_{j=1}^\ell g_j\big(\omega_{s_j}-\omega_0-L_{s_j}(\omega,\widetilde{\mu}^N)\big)\bigg\rangle\!\bigg)\big(1-\langle\widetilde{\mu}^N, m_{T'}\rangle\big)\bigg],
\end{equation}
where
\begin{eqnarray}
&& L_t(\omega,\mu):=((\omega_t+1)^+\wedge 1)\,C\log\big(1-\langle \mu,m_t\rangle\big), \\
&& \widetilde{G}(y) :=\bigg(\!G(y) - G\bigg(\!\EE\bigg[\prod_{j=1}^\ell (\alpha\,s_j+\sigma\,\overline{B}_{s_j})\bigg]\bigg)\!\bigg)^2,
\end{eqnarray}
and we use the convention that the expression inside the expectation in the definition of $Q^N$ is $0$ whenever $\langle \widetilde{\mu}^N, m_{T'}\rangle=1$. Note 
\begin{equation}
Q^N = \EE\bigg[\widetilde{G}\bigg(\!\frac{1}{N}\,\sum_{i=1}^N \,\prod_{j=1}^\ell (\alpha\,s_j+\sigma\,B^i_{s_j})\!\bigg)\big(1-\langle \tilde{\mu}^N, m_{T'}\rangle\big)\bigg],
\end{equation}
so that $\lim_{N\to\infty} Q^N=0$ by the strong law of large numbers. The last ingredient in the proof of Theorem \ref{main2} is the following lemma, which is the analogue of \cite[Lemma 5.10]{Delarue2}.

\begin{lemma}\label{le:cont.main}
The functional
\begin{equation}
\mathcal{P}(\mathcal{D}([0,T+1]))\ni\mu \mapsto \widetilde{G}\bigg(\!\bigg\langle \mu,\,
\prod_{j=1}^\ell g_j\big(\omega_{s_j}-\omega_0-L_{s_j}(\omega,\mu)\big)\bigg\rangle\!\bigg)\big(1-\langle\mu, m_{T'}\rangle\big)
\end{equation}
is continuous at $\Pi_{\infty}$-almost every $\mu$.
\end{lemma}

\smallskip

\noindent\textbf{Proof.} Lemma \ref{le:cont.mum} implies that the mappings $\mu\mapsto\langle\mu,m_{s_j}\rangle$, $j=1,\,2,\,\ldots,\,\ell$ and $\mu\mapsto\langle\mu,m_{T'}\rangle$ are continuous at $\Pi_{\infty}$-almost every $\mu$. Pick such a $\mu$ and $\mu_n$, $n\in\nn$ converging to $\mu$. If $\langle \mu, m_{s_\ell}\rangle = 1$, then $\langle \mu, m_{T'}\rangle = 1$, and the value of the functional at $\mu$ is $0$. At the same time, $\lim_{n\to\infty} \langle \mu_n, m_{T'}\rangle=1$, and the boundedness of $\widetilde{G}$ implies that the values of the functional at $\mu_n$, $n\in\nn$ converge to $0$, yielding the desired continuity.

\medskip

If $\langle \mu, m_{s_\ell}\rangle < 1$, then for all $n\in\nn$ sufficiently large $\langle \mu_n, m_{s_\ell}\rangle$ is bounded away from $1$ by a constant (recall the continuity of $\langle\cdot,m_{s_\ell}\rangle$ at $\mu$). In particular, no discontinuity can arise from the logarithm in the definition of $L_t(\omega,\mu)$. Consequently, we can repeat the proof of \cite[Lemma 5.10]{Delarue2} to conclude that the values of the functional at $\mu_n$, $n\in\nn$ converge to its value at $\mu$. \ep

\medskip

Lemma \ref{le:cont.main} gives
\begin{equation}
\int_{\mathcal{P}(\mathcal{D}([0,T+1]))} 
\widetilde{G}\bigg(\!\bigg\langle \mu,\,
\prod_{j=1}^\ell g_j\big(\omega_{s_j}-\omega_0-L_{s_j}(\omega,\mu)\big)\bigg\rangle\!\bigg)\big(1-\langle\mu, m_{T'}\rangle\big)\,\Pi_\infty(\mathrm{d}\mu)\!=\!\lim_{N\rightarrow\infty} Q^N \!=\! 0,
\end{equation}
and, hence,
\begin{equation}\label{G_conc}
G\bigg(\!\bigg\langle \mu,\,
\prod_{j=1}^\ell g_j\big(\omega_{s_j}-\omega_0-L_{s_j}(\omega,\mu)\big)\bigg\rangle\!\bigg) = 
G\bigg(\!\EE\bigg[\prod_{j=1}^\ell (\alpha\,s_j+\sigma\,\overline{B}_{s_j})\bigg]\bigg)
\end{equation}
for $\Pi_{\infty}$-almost every $\mu$ with $\langle \mu, m_{T'}\rangle < 1$. The standard arguments in the proof of \cite[Lemma 5.4]{Delarue2} allow to deduce from \eqref{G_conc} that the process $\omega_t-\omega_0-L_t(\omega,\mu)$, $t\in[0,T']$ is a Brownian motion with drift and $\omega_0\stackrel{d}{=}\widetilde{Y}^1_0$, under $\Pi_{\infty}$-almost every $\mu$ with $\langle \mu, m_{T'}\rangle < 1$. Since the set of possible $T'$ is countable and dense in $[0,T]$, we conclude that the canonical process satisfies the condition \eqref{eq:phys_sol.1} in Definition \ref{def_phys_sol.new} of a physical solution, under $\Pi_{\infty}$-almost every $\mu$.

\medskip

To see the condition \eqref{eq:phys_sol.2} in Definition \ref{def_phys_sol.new} we cast the estimate \eqref{eq:jump.necessary} of Lemma \ref{le:jump.necessary} as
\begin{equation}
\begin{split}
& \;\PP\bigg(\langle \widetilde{\mu}^{N},m_{t-} \rangle\leq 1-r,\;
\forall\iota \leq \bigg(\frac{\langle \widetilde{\mu}^{N}, m_{t+s} \rangle-\langle \widetilde{\mu}^{N}, m_{t-} \rangle}{1-\langle \widetilde{\mu}^{N}, m_{t-} \rangle}-2s^{1/4}\bigg)^+:\\ 
& \qquad\qquad\qquad\qquad\qquad\;
\frac{\widetilde{\mu}^N(m_{t-}=0,\,\omega_{t-}+C\log(1-\iota)\leq 2s^{1/4})}
{1-\langle \widetilde{\mu}^{N}, m_{t-} \rangle} \geq \frac{\iota}{1+s^{1/4}}\bigg) \\
& \geq\PP\big(\langle\widetilde{\mu}^N,m_{t-}\rangle\leq 1-r\big)-C_0s.
\end{split}
\end{equation}
By following the last part of the proof of \cite[Theorem 4.4]{Delarue2} we obtain from this for all $r\in(0,1)$ sufficiently small and all $t\in J\cap[0,T]$: 
\begin{equation}
\begin{split}
& \;\Pi_\infty\bigg(\!\langle \mu, m_t \rangle\leq 1\!-\!r/2,\,\forall\iota<\frac{\langle \mu,m_t\!-\!m_{t-}\rangle}{1\!-\!\langle\mu,m_{t-}\rangle}:\,
\frac{\mu(m_{t-}\!=\!0,\,\omega_{t-}\!+\!C\log(1\!-\!\iota)\!\leq\!0)}
{1\!-\!\langle \mu, m_{t-} \rangle} \geq \iota\!\bigg) \\ 
& \geq \Pi_\infty\big(\langle \mu, m_t \rangle \leq 1 - r \big),
\end{split}
\end{equation}
which in the limit $r\downarrow0$ yields
\begin{equation}
\begin{split}
& \;\Pi_\infty\bigg(\!\langle \mu, m_t\rangle < 1,\,\forall\iota < \frac{\langle \mu,m_t\!-\!m_{t-}\rangle}{1\!-\!\langle\mu,m_{t-} \rangle}:\,
\frac{\mu(m_{t-}\!=\!0,\,\omega_{t-}\!+\!C\log(1\!-\!\iota)\!\leq\!0)}
{1\!-\!\langle \mu, m_{t-} \rangle} \geq \iota \!\bigg) \\
& \geq \Pi_\infty\big(\langle \mu, m_t \rangle < 1 \big).
\end{split}
\end{equation}
Consequently, 
\begin{equation}
\forall\iota<\frac{\langle\mu,m_t-m_{t-}\rangle}{1-\langle \mu, m_{t-} \rangle}:\;
\frac{\mu(m_{t-}=0,\,\omega_{t-}+C\log(1-\iota)\leq 0)}
{1-\langle\mu,m_{t-}\rangle}\geq\iota
\end{equation}
for all $t\in J\cap [0,T]$ and $\Pi_\infty$-almost every $\mu$ with 
\begin{equation}
\overline{\tau}^0(\mu):=\inf\{s\in[0,T]:\langle \mu, m_s \rangle=1\}>t.
\end{equation} 
Since the set $J$ is countable, the latter statement holds for $\Pi_\infty$-almost every $\mu$ and all $t\in J\cap [0,T]\cap [0,\overline{\tau}^0(\mu))$. This implies that the canonical process satisfies the condition \eqref{eq:phys_sol.2} in Definition \ref{def_phys_sol.new} of a physical solution, under $\Pi_\infty$-almost every $\mu$. The proof of Theorem \ref{main2} is finished.

\subsection{Proof of Lemma \ref{le:cont.mum}}
\label{subse:cont.mum}

Let us fix an arbitrary $T_1\in(0,\infty)$ and, as before, denote by $\mathcal{D}([0,T_1])$ the space of real-valued c\`adl\`ag functions on $[0,T_1]$ that are continuous at $T_1$, endowed with the Skorokhod M1 topology. The following two lemmas are the analogues of \cite[Lemma 5.6]{Delarue2} and \cite[Proposition 5.8]{Delarue2}, respectively, and we omit their proofs, since they constitute very minor modifications of the proofs in \cite{Delarue2}.

\begin{lemma}\label{le.Delarue2.lemma5.6}
Consider any $\omega\in\mathcal{D}([0,T_1])$ satisfying the crossing property
\begin{equation}
\forall\,s>0:\;\;\tau<T_1\,\Longrightarrow\, \inf_{t\in[\tau,(\tau+s)\wedge T_1]}\,(\omega_t - \omega_\tau) < 0,
\end{equation}
where $\tau:= (\inf\{t\in[0,T_1]:\,\omega_t\leq 0\})\wedge T_1$. Then, for any $\omega^n$, $n\in\nn$ converging to $\omega$ in $\mathcal{D}([0,T_1])$ there exists a countable set $J\subset [0,T_1]$ such that
\begin{equation}
\lim_{n\to\infty} m^n_t:=\lim_{n\to\infty} \bone_{\{\inf_{s\in[0,t]} \omega^n_s \leq 0\}}
=\bone_{\{\inf_{s\in[0,t]} \omega_s \leq 0\}}=:m_t,\;\; t\in[0,T_1]\setminus J
\end{equation}
and all points $t$ of continuity of $\omega$ satisfying $\inf_{s\in[0,t]} \omega_s \neq 0$ are contained in $[0,T_1]\setminus J$. 
\end{lemma}

\begin{lemma}\label{le.Delarue2.prop5.8}
Consider any $\mu\in \mathcal{P}(\mathcal{D}([0,T_1]))$ satisfying
\begin{equation}
\forall\,s>0:\;\;\mu\Big(\tau<T_1,\,\inf_{t\in[\tau,(\tau+s)\wedge T_1]} \,(\omega_t - \omega_{\tau}) = 0 \Big) = 0,
\end{equation}
where $\tau:= (\inf\{t\in[0,T_1]:\,\omega_t\leq 0\})\wedge T_1$. Then, for any $\mu_n$, $n\in\nn$ converging weakly to $\mu$ we have
\begin{equation}
\lim_{n\rightarrow\infty}\,\langle \mu_n, m_t \rangle = \langle \mu, m_t \rangle
\end{equation}
for all points $t$ of continuity of the mapping $t\mapsto \langle \mu, m_t \rangle$.
\end{lemma}

The next lemma shows that the canonical process satisfies the crossing property under $\Pi_\infty$-almost every $\mu$. It is the analogue of \cite[Lemma 5.9]{Delarue2}.

\begin{lemma}\label{le:Delarue2.lemma5.9}
For $\Pi_\infty$-almost every $\mu$ it holds
\begin{equation}
\forall\,s>0:\;\;\mu\Big(\tau<T+1,\,\inf_{t\in[\tau,(\tau+s)\wedge(T+1)]}\, (\omega_t - \omega_\tau)=0\Big)=0,
\end{equation}
where $\tau:=(\inf\{t\in[0,T+1]\,:\,\omega_t\leq 0\})\wedge (T+1)$.
\end{lemma}

\noindent The proof of Lemma \ref{le:Delarue2.lemma5.9} is essentially the same as the proof of \cite[Lemma 5.9]{Delarue2}, with Step 2 therein allowing for a simplification, since the drift coefficient $\alpha$ is constant in the present setting. It is also worth mentioning the typo in the seventh displayed equation in the proof of \cite[Lemma 5.9]{Delarue2}: ``$\max(\zeta^N_{r^N_t-},\zeta^N_{r^N_t}) - \min(\zeta^N_{r^N_s-},\zeta^N_{r^N_s})$" should be replaced by ``$\min(\zeta^N_{r^N_t-},\zeta^N_{r^N_t}) - \max(\zeta^N_{r^N_s-},\zeta^N_{r^N_s})$".

\medskip

Lastly, we observe that Lemma \ref{le:cont.mum} is a direct consequence of Lemmas \ref{le.Delarue2.prop5.8} and \ref{le:Delarue2.lemma5.9}.




\bigskip\bigskip\bigskip

\bibliographystyle{amsalpha}
\bibliography{SystemicRisk}

\bigskip\bigskip\bigskip

\end{document}